\def\enroute{y}
\def\draft{n}

\documentclass[12pt]{amsart}
\usepackage{amssymb,epic,eepic,amscd,xr}
\theoremstyle{plain}

\newtheorem{theorem}{Theorem}
\newtheorem{proposition}{Proposition}[section]
\newtheorem{fact}[proposition]{Fact}

\newtheorem{corollary}[proposition]{Corollary}
\newtheorem{claim}[proposition]{Claim}
\newtheorem{Truth}[proposition]{Truth}
\newtheorem{AlmostTruth}[proposition]{Almost Truth}

\newtheorem{slogan}{Slogan}

\theoremstyle{definition}
\newtheorem{definition}[proposition]{Definition}

\newtheorem{question}[proposition]{Question}

\theoremstyle{remark}

\newtheorem{exercise}[proposition]{Exercise}

\def\printname#1{
	\if\draft y
		\smash{\makebox[0pt]{\hspace{-0.5in}
			\raisebox{8pt}{\tt\tiny #1}}}
	\fi
}

\newcommand{\mathmode}[1]{$#1$}

\newlength{\standardunitlength}
\catcode`\@=11
\long\def\@makecaption#1#2{%
    \vskip 10pt
    \setbox\@tempboxa\hbox{
      \small\sf{\bfcaptionfont #1. }\ignorespaces #2}%
    \ifdim \wd\@tempboxa >\captionwidth {%
        \rightskip=\@captionmargin\leftskip=\@captionmargin
        \unhbox\@tempboxa\par}%
      \else
        \hbox to\hsize{\hfil\box\@tempboxa\hfil}%
    \fi}
\font\bfcaptionfont=cmssbx10 scaled \magstephalf
\newdimen\@captionmargin\@captionmargin=2\parindent
\newdimen\captionwidth\captionwidth=\hsize
\catcode`\@=12

\newcommand{\Ad}{\operatorname{Ad}}

\newcommand{\gr}{\operatorname{gr}}

\newcommand{\tr}{\operatorname{tr}}

\def\lbl#1{\label{#1}\printname{#1}}

\newlength{\globalparindent}
\newcommand{\udot}{{\mathaccent\cdot\cup}}

\advance \topmargin by -\headheight
\advance \topmargin by -\headsep
\evensidemargin \oddsidemargin
\newcommand{\Aa}{\text{\sl\AA}}
\newcommand{\Arhus}{\AA rhus}

\newcommand{\BCH}{\text{BCH}}

\newcommand{\Gint}{\int^{FG}}
\newcommand{\LMO}{\text{LMO}}

\newcommand{\Tg}{{\mathcal T}_\frakg}

\newcommand{\Zcheck}{\check{Z}}
\newcommand{\calA}{\mathcal A}
\newcommand{\calB}{\mathcal B}
\newcommand{\calD}{\mathcal D}

\newcommand{\calN}{\mathcal N}
\newcommand{\frakg}{\mathfrak g}

\newcommand{\Answer}{\par\noindent {\em Answer: }}

\newcommand\ifenroute[2]{\if\enroute y #1 \else #2 \fi}

\newcommand\chordn[2]{\,^{#1}\!\!\frown^{#2}}
\newcommand\chordu[2]{\,_{#1}\!\!\smile_{#2}}

\begin{document}

\title[The \AA{}rhus integral I: Introductions]
  {
  The \AA{}rhus integral of rational homology 3-spheres I: A highly non
  trivial flat connection on $S^3$.
}

\author[Bar-Natan]{Dror~Bar-Natan}
\address{Institute of Mathematics\\
        The Hebrew University\\
        Giv'at-Ram, Jerusalem 91904\\
        Israel}
\email{drorbn@math.huji.ac.il}

\author[Garoufalidis]{Stavros~Garoufalidis}
\address{Department of Mathematics\\
        Harvard University\\
        Cambridge MA 02138\\
        USA}
\email{stavros@math.harvard.edu}

\author[Rozansky]{Lev~Rozansky}
\address{Department of Mathematics, Statistics, and Computer Science\\
        University of Illinois at Chicago\\
        Chicago IL 60607-7045\\
        USA}
\curraddr{Department of Mathematics\\
  Yale University\\
  10 Hillhouse Avenue\\
  P.O. Box 208283 \\
  New Haven, CT 06520-8283\\
  USA
}
\email{rozansky@math.yale.edu}

\author[Thurston]{Dylan~P.~Thurston}
\address{Department of Mathematics\\
        University of California at Berkeley\\
        Berkeley CA 94720-3840\\
        USA}
\email{dpt@math.berkeley.edu}

\thanks{This paper is available electronically at
  {\tt http://www.ma.huji.ac.il/\~{}drorbn}, at \newline
  {\tt http://jacobi.math.brown.edu/\~{}stavrosg}, and at {\tt
    http://xxx.lanl.gov/abs/q-alg/9706004}.
}

\dedicatory{To appear in {\em Selecta Mathematica.}}
\date{This edition: Feb.~15,~1999; \ \ First edition: June~4, 1997.}

\begin{abstract}
Path integrals don't really exist, but it is very useful to dream that
they do exist, and figure out the consequences. Apart from describing
much of the physical world as we now know it, these dreams also lead to
some highly non-trivial mathematical theorems and theories. We argue that
even though non-trivial flat connections on $S^3$ don't really exist,
it is beneficial to dream that one exists (and, in fact, that it comes
from the non-existent Chern-Simons path integral). Dreaming the right
way, we are led to a rigorous construction of a universal finite-type
invariant of rational homology spheres.  We show that this invariant
is equal (up to a normalization) to the $\LMO$ (Le-Murakami-Ohtsuki,
\cite{LeMurakamiOhtsuki:Universal}) invariant and that it recovers the
Rozansky and Ohtsuki invariants.

This is part I of a 4-part series, containing the introductions and
answers to some frequently asked questions. Theorems are stated but not
proved in this part, and it can be viewed as a ``research announcement''.
Part II of this series is titled ``Invariance and Universality'', part
III ``The Relation with the Le-Murakami-Ohtsuki Invariant'', and part
IV ``The Relation with the Rozansky and Ohtsuki Invariants''.
\end{abstract}

\maketitle

\tableofcontents

{\it \ifenroute{
  This series has two introductions. If your flight is about to land and
  you have only little time for reading, read only the first one; it's
  more fun.  We are confident that after that, you will not be able to
  resist the temptation to read the other in the taxi on the way home.

  On the other hand, if you have a perverse aversion to philosophy, the
  second introduction and the rest of the series are fully rigorous and
  can be read independently of the first introduction.
}{
  This series has two introductions. The first is philosophical
  and non-rigorous. We recommend reading it first. The second introduction
  and the rest of the series are fully rigorous and can be read
  independently of the first introduction.
}}

\section{Introduction for philosophy majors} \label{sec:intro1}

\subsection{What if there were?} \lbl{subsec:WhatIf}

Suppose there was some highly-non-trivial flat connection $A$ on
$S^3$. Well, of course there are no non-trivial flat connections on $S^3$;
in recent years there has been no significant debate over that. But we wrote
``highly-non-trivial flat connection''; meaning that the object that
we're going to talk about is neither a connection nor it is flat (and
possibly, it doesn't even exist). But we will see that it is beneficial
to assume that such an $A$ does exist, and that it has some fixed good
properties. We will make some deductions and guess some formulas, and
later on we will prove that they work, with no reference to $A$.

We claim that given a well-behaved highly-non-trivial flat connection
$A$ on $S^3$, we can construct (under some mild conditions) a link
invariant\footnote{
  \lbl{foot:Aa}
  \Aa{} is the last letter of the Danish alphabet.  It is pronounced like
  the ``O'' in ``Oval''.  In typographically deprived environments, it
  may be spelled Aa.
}
$\Aa_A$ that respects the Kirby moves, and hence an invariant
of 3-manifolds. Just constructing a link invariant is easy; simply
consider the holonomy $h_A(L)$ of the connection $A$ along a link
$L$. The invariance under small deformations (that do not pass through
self-intersections) of an embedding of $L$ follows from the flatness of
$A$. The non-invariance under full homotopy comes from the fact that $A$ is
highly-non-trivial. If $h_A(L)$ were invariant under homotopy, it would
have been rather dull, for all links of a fixed number of components
are homotopic.

The construction of $\Aa_A(L)$ from $h_A(L)$ is extremely simple, and
can be summarized by a single catchy slogan\ifenroute{(printed in large
letters, in case there's turbulence):}{:}

\begin{slogan} \lbl{IntegrateHolonomies}
\[ \fbox{\text{\Large\it Integrate the Holonomies}} \]
\end{slogan}

Let's try to make sense out of that. Suppose $A$ is a $\frakg$-connection
for some metrized Lie-algebra $\frakg$. The holonomy of $A$ along
a single path is roughly the product of the values of $A$ seen
along the path. Namely, it is a certain long product of elements of
$\frakg$. So it is naturally valued in $\hat{U}(\frakg)$, a certain
completion of the universal enveloping algebra of $\frakg$. By the
Poincare-Birkhoff-Witt theorem, $\hat{U}(\frakg)$ is isomorphic to
$\hat{S}(\frakg)$, a completion of the symmetric algebra of $\frakg$, via
the symmetrization map $\beta:\hat{S}(\frakg)\to\hat{U}(\frakg)$ (see
e.g.~\cite[paragraph~2.4.10]{Dixmier:EnvelopingAlgebras}). The algebra
$\hat{S}(\frakg)$ is the algebra of power series on $\frakg^\star$,
and those power series that are convergent are called functions and can
sometimes be integrated. It is in this sense that one should interpret
slogan~\ref{IntegrateHolonomies}, only that in the case of an $X$-marked
link (a link whose components are in a bijective correspondence with some
$n$-element set of labels $X=\{x_i\}_{i=1}^n$, with a base point on
each component, to make the holonomies well defined) the holonomies are
in $\hat{U}(\frakg)^{\otimes n}$, and thus the integration is over $n$
copies of $\frakg^\star$, with one integration variable (also denoted
$x_i$) for each component of the link:

\begin{definition} \lbl{ArhusDefinition} The \Arhus\footnote{
\Arhus{} ($\hat{O}R{\mathbf '}h\overline{oo}s'$), n, also Aarhus.  A city
of central Denmark, a commercial and industrial center on \Arhus{} Bay, an
arm of the Kattegat. First mentioned in the mid-tenth century, \Arhus{}
is one of the country's oldest cities. Population, 250,404. ({\em The
American Heritage Dictionary of the English Language, Third Edition},
copyright \copyright 1992 by Houghton Mifflin Company.)}
integral $\Aa_A$ for the non-trivial
flat connection $A$ is the integral
\[ \Aa_A(L)=
  \calN \int_{\frakg^\star\oplus\ldots\oplus\frakg^\star}
  h_A(L)(x_1,x_2,\ldots,x_n)
  \,dx_1\,dx_2\ldots dx_n,
\]
where $h_A(L)(x_1,x_2,\ldots,x_n)$ is the holonomy $h_A(L)$ regarded
as a function of
$(x_1,\ldots,x_n)\in\frakg^\star\oplus\ldots\oplus\frakg^\star$, the
symbol $dx$ denotes the measure on $\frakg^\star$ induced by its metric,
and $\calN$ denotes additional normalizations which are of interest only
for math majors.
\end{definition}

Let us get to the main point \ifenroute{before, in preparation for
landing, you will be told to make sure that your tray-table is in its
fully upright and locked position. Reading will be harder then.}{as
quickly as possible.} Why is $\Aa_A$ invariant under the second and
more difficult Kirby move~\cite{Kirby:Calculus}? Not worrying too much
about the important issue of framing (we can do it later, when we're
on the ground), the second Kirby move is the operation of `sliding' one
component of the link along a neighboring one. See figure~\ref{KirbyII}.

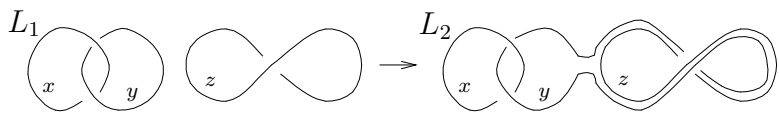
\begin{figure}[htpb]
\[ \printname{KirbyII}
	\setlength{\unitlength}{0.75\standardunitlength}
	\begin{array}{c}  \hspace{-1.7mm}
        	\raisebox{-8pt}{\begingroup\makeatletter\ifx\SetFigFont\undefined
\def\x#1#2#3#4#5#6#7\relax{\def\x{#1#2#3#4#5#6}}%
\expandafter\x\fmtname xxxxxx\relax \def\y{splain}%
\ifx\x\y   
\gdef\SetFigFont#1#2#3{%
  \ifnum #1<17\tiny\else \ifnum #1<20\small\else
  \ifnum #1<24\normalsize\else \ifnum #1<29\large\else
  \ifnum #1<34\Large\else \ifnum #1<41\LARGE\else
     \huge\fi\fi\fi\fi\fi\fi
  \csname #3\endcsname}%
\else
\gdef\SetFigFont#1#2#3{\begingroup
  \count@#1\relax \ifnum 25<\count@\count@25\fi
  \def\x{\endgroup\@setsize\SetFigFont{#2pt}}%
  \expandafter\x
    \csname \romannumeral\the\count@ pt\expandafter\endcsname
    \csname @\romannumeral\the\count@ pt\endcsname
  \csname #3\endcsname}%
\fi
\fi\endgroup
{\renewcommand{\dashlinestretch}{30}
\begin{picture}(6406,865)(0,-10)
\path(3097,392)(3397,392)
\path(3277.000,362.000)(3397.000,392.000)(3277.000,422.000)
\path(697,542)	(673.470,501.490)
	(654.424,465.836)
	(628.904,406.464)
	(618.682,358.610)
	(622.000,317.000)

\path(622,317)	(637.171,267.719)
	(659.699,217.347)
	(689.097,168.112)
	(724.874,122.244)
	(766.542,81.970)
	(813.611,49.519)
	(865.594,27.119)
	(922.000,17.000)

\path(922,17)	(987.928,19.472)
	(1052.134,33.650)
	(1112.757,58.686)
	(1167.936,93.729)
	(1215.811,137.931)
	(1254.520,190.443)
	(1282.204,250.415)
	(1297.000,317.000)

\path(1297,317)	(1296.200,384.584)
	(1279.981,447.654)
	(1250.809,505.509)
	(1211.151,557.446)
	(1163.475,602.766)
	(1110.246,640.765)
	(1053.932,670.744)
	(997.000,692.000)

\path(997,692)	(954.226,697.450)
	(905.984,687.939)
	(847.000,661.708)
	(772.000,617.000)

\path(772,167)	(793.571,208.705)
	(811.218,245.211)
	(835.619,305.266)
	(846.960,352.438)
	(847.000,392.000)

\path(847,392)	(827.499,453.439)
	(790.131,517.550)
	(743.698,575.136)
	(697.000,617.000)

\path(697,617)	(635.379,654.208)
	(597.963,671.940)
	(557.882,687.095)
	(516.442,698.175)
	(474.945,703.684)
	(434.696,702.125)
	(397.000,692.000)

\path(397,692)	(348.071,665.503)
	(303.573,629.428)
	(264.261,585.653)
	(230.894,536.056)
	(204.227,482.515)
	(185.018,426.906)
	(174.023,371.109)
	(172.000,317.000)

\path(172,317)	(178.582,271.074)
	(193.084,224.563)
	(214.623,178.963)
	(242.316,135.770)
	(312.640,62.592)
	(353.507,35.600)
	(397.000,17.000)

\path(397,17)	(440.100,10.914)
	(488.451,20.214)
	(547.327,46.656)
	(622.000,92.000)

\path(2122,467)	(2072.998,522.548)
	(2029.039,568.943)
	(1989.246,606.844)
	(1952.740,636.909)
	(1886.071,676.170)
	(1822.000,692.000)

\path(1822,692)	(1781.624,691.861)
	(1739.129,684.935)
	(1695.987,671.958)
	(1653.666,653.664)
	(1613.637,630.786)
	(1577.370,604.060)
	(1522.000,542.000)

\path(1522,542)	(1494.504,472.878)
	(1487.629,433.174)
	(1485.336,392.000)
	(1487.627,350.826)
	(1494.501,311.122)
	(1522.000,242.000)

\path(1522,242)	(1577.182,179.694)
	(1613.281,152.749)
	(1653.156,129.673)
	(1695.387,111.262)
	(1738.553,98.315)
	(1781.231,91.629)
	(1822.000,92.000)

\path(1822,92)	(1877.607,106.804)
	(1931.513,136.486)
	(1983.321,176.849)
	(2032.634,223.696)
	(2079.055,272.832)
	(2122.188,320.059)
	(2197.000,392.000)

\path(2197,392)	(2238.717,424.983)
	(2287.976,466.761)
	(2343.083,513.583)
	(2402.342,561.695)
	(2464.059,607.346)
	(2526.538,646.782)
	(2588.083,676.251)
	(2647.000,692.000)

\path(2647,692)	(2703.554,690.790)
	(2766.535,675.826)
	(2825.999,650.199)
	(2872.000,617.000)

\path(2872,617)	(2903.943,570.158)
	(2927.477,509.293)
	(2942.023,446.031)
	(2947.000,392.000)

\path(2947,392)	(2942.023,337.969)
	(2927.477,274.707)
	(2903.943,213.842)
	(2872.000,167.000)

\path(2872,167)	(2825.840,134.044)
	(2766.111,108.822)
	(2703.077,93.940)
	(2647.000,92.000)

\path(2647,92)	(2576.950,109.023)
	(2538.746,125.693)
	(2496.959,148.681)
	(2450.489,178.647)
	(2398.239,216.249)
	(2339.108,262.147)
	(2272.000,317.000)

\path(4223,167)	(4244.571,208.705)
	(4262.218,245.211)
	(4286.619,305.266)
	(4297.960,352.438)
	(4298.000,392.000)

\path(4298,392)	(4278.500,453.439)
	(4241.135,517.550)
	(4194.703,575.136)
	(4148.000,617.000)

\path(4148,617)	(4086.379,654.208)
	(4048.963,671.940)
	(4008.883,687.095)
	(3967.442,698.175)
	(3925.945,703.684)
	(3885.696,702.125)
	(3848.000,692.000)

\path(3848,692)	(3799.071,665.503)
	(3754.573,629.428)
	(3715.261,585.653)
	(3681.894,536.056)
	(3655.227,482.515)
	(3636.018,426.906)
	(3625.023,371.109)
	(3623.000,317.000)

\path(3623,317)	(3629.582,271.074)
	(3644.084,224.563)
	(3665.623,178.963)
	(3693.316,135.770)
	(3763.640,62.592)
	(3804.507,35.600)
	(3848.000,17.000)

\path(3848,17)	(3891.100,10.914)
	(3939.451,20.214)
	(3998.327,46.656)
	(4073.000,92.000)

\path(5573,467)	(5523.998,522.548)
	(5480.039,568.943)
	(5440.246,606.844)
	(5403.740,636.909)
	(5337.071,676.170)
	(5273.000,692.000)

\path(5273,692)	(5232.624,691.861)
	(5190.129,684.935)
	(5146.987,671.958)
	(5104.666,653.664)
	(5064.637,630.786)
	(5028.370,604.060)
	(4973.000,542.000)

\path(4973,542)	(4945.504,472.878)
	(4938.629,433.174)
	(4936.336,392.000)
	(4938.627,350.826)
	(4945.501,311.122)
	(4973.000,242.000)

\path(4973,242)	(5028.182,179.694)
	(5064.281,152.749)
	(5104.156,129.673)
	(5146.387,111.262)
	(5189.553,98.315)
	(5232.231,91.629)
	(5273.000,92.000)

\path(5273,92)	(5328.607,106.804)
	(5382.513,136.486)
	(5434.321,176.849)
	(5483.634,223.696)
	(5530.055,272.832)
	(5573.188,320.059)
	(5648.000,392.000)

\path(5648,392)	(5689.717,424.983)
	(5738.976,466.761)
	(5794.083,513.583)
	(5853.343,561.695)
	(5915.059,607.346)
	(5977.538,646.782)
	(6039.083,676.251)
	(6098.000,692.000)

\path(6098,692)	(6154.554,690.790)
	(6217.535,675.826)
	(6276.999,650.199)
	(6323.000,617.000)

\path(6323,617)	(6354.943,570.158)
	(6378.477,509.293)
	(6393.023,446.031)
	(6398.000,392.000)

\path(6398,392)	(6393.023,337.969)
	(6378.477,274.707)
	(6354.943,213.842)
	(6323.000,167.000)

\path(6323,167)	(6276.840,134.044)
	(6217.111,108.822)
	(6154.077,93.940)
	(6098.000,92.000)

\path(6098,92)	(6027.950,109.023)
	(5989.746,125.693)
	(5947.959,148.681)
	(5901.489,178.647)
	(5849.239,216.249)
	(5790.108,262.147)
	(5723.000,317.000)

\path(4747,467)	(4711.461,518.864)
	(4679.358,562.626)
	(4622.826,628.486)
	(4572.131,669.853)
	(4522.000,692.000)

\path(4522,692)	(4469.151,698.609)
	(4406.222,689.484)
	(4368.781,678.472)
	(4326.182,662.867)
	(4277.549,642.449)
	(4222.000,617.000)

\path(4148,542)	(4124.470,501.490)
	(4105.424,465.836)
	(4079.904,406.464)
	(4069.682,358.610)
	(4073.000,317.000)

\path(4073,317)	(4089.311,268.601)
	(4114.609,220.372)
	(4147.427,173.784)
	(4186.299,130.310)
	(4229.758,91.422)
	(4276.340,58.593)
	(4324.575,33.295)
	(4373.000,17.000)

\path(4373,17)	(4431.290,16.459)
	(4493.625,34.393)
	(4551.647,62.380)
	(4597.000,92.000)

\path(4597,92)	(4661.209,164.986)
	(4700.012,228.478)
	(4722.637,269.281)
	(4748.000,317.000)

\path(5767,347)	(5829.436,296.293)
	(5876.228,260.058)
	(5939.000,219.000)

\path(5939,219)	(5972.340,203.274)
	(6015.064,186.231)
	(6059.255,172.573)
	(6097.000,167.000)

\path(6097,167)	(6145.553,170.388)
	(6202.515,181.574)
	(6257.220,202.222)
	(6299.000,234.000)

\path(6299,234)	(6316.570,269.215)
	(6323.944,312.640)
	(6324.595,356.245)
	(6322.000,392.000)

\path(6322,392)	(6315.455,440.152)
	(6302.451,497.805)
	(6280.472,553.305)
	(6247.000,595.000)

\path(6247,595)	(6196.420,618.872)
	(6135.220,628.204)
	(6073.660,626.434)
	(6022.000,617.000)

\path(6022,617)	(5984.285,601.944)
	(5943.865,578.452)
	(5902.125,548.874)
	(5860.446,515.562)
	(5820.212,480.868)
	(5782.804,447.142)
	(5722.000,392.000)

\path(5722,392)	(5671.579,343.226)
	(5610.111,279.020)
	(5548.338,215.054)
	(5497.000,167.000)

\path(5497,167)	(5452.204,128.322)
	(5394.944,80.998)
	(5332.462,39.175)
	(5272.000,17.000)

\path(5272,17)	(5215.585,20.699)
	(5152.836,39.811)
	(5093.419,66.268)
	(5047.000,92.000)

\path(5047,92)	(5008.855,120.610)
	(4965.599,160.505)
	(4925.544,203.648)
	(4897.000,242.000)

\path(4897,242)	(4889.738,279.785)
	(4874.000,317.000)

\path(4874,317)	(4826.944,331.250)
	(4792.571,327.688)
	(4747.000,317.000)

\path(4747,467)	(4796.122,458.597)
	(4832.886,457.005)
	(4882.000,475.000)

\path(4882,475)	(4897.000,535.000)

\path(4897,535)	(4925.029,575.285)
	(4964.894,620.411)
	(5008.312,662.081)
	(5047.000,692.000)

\path(5047,692)	(5093.398,718.044)
	(5152.818,744.466)
	(5215.578,763.405)
	(5272.000,767.000)

\path(5272,767)	(5332.334,744.576)
	(5394.602,702.339)
	(5451.820,654.932)
	(5497.000,617.000)

\path(5497,617)	(5542.076,579.202)
	(5576.944,547.964)
	(5624.000,505.000)

\put(292,167){\makebox(0,0)[lb]{\smash{{\mathmode{\scriptstyle x}}}}}
\put(989,121){\makebox(0,0)[lb]{\smash{{\mathmode{\scriptstyle y}}}}}
\put(3750,167){\makebox(0,0)[lb]{\smash{{\mathmode{\scriptstyle x}}}}}
\put(4417,121){\makebox(0,0)[lb]{\smash{{\mathmode{\scriptstyle y}}}}}
\put(5076,227){\makebox(0,0)[lb]{\smash{{\mathmode{\scriptstyle z}}}}}
\put(1642,219){\makebox(0,0)[lb]{\smash{{\mathmode{\scriptstyle z}}}}}
\put(0,670){\makebox(0,0)[lb]{\smash{{\mathmode{L_1}}}}}
\put(3419,647){\makebox(0,0)[lb]{\smash{{\mathmode{L_2}}}}}
\end{picture}
} }
        	\hspace{-1.9mm}
	\end{array}
 \]
\caption{The second Kirby move.}
\lbl{KirbyII}
\end{figure}

Let us analyze the behavior of $\Aa_A$ under the second Kirby move. The
only difference is that the holonomy along the component labeled $y$
changes, and (with an appropriate choice of basepoints)
the change is rather simple: it gets multiplied by the
holonomy around the component labeled $z$. In formulas valid in
$\hat{U}(\frakg)^{\otimes 3}$, this amounts to saying
\begin{equation} \lbl{Uglevel}
  h_A(L_2)=(1\otimes \times_U\otimes 1)(1\otimes 1 \otimes\Delta)h_A(L_1),
\end{equation}
where $\Delta:\hat{U}(\frakg)\to\hat{U}(\frakg)\otimes\hat{U}(\frakg)$
is the co-product, which in~\eqref{Uglevel} takes the
group-like holonomy of $A$ along $z$ and ``doubles'' it, and
$\times_U:\hat{U}(\frakg)\otimes\hat{U}(\frakg)\to\hat{U}(\frakg)$ is
the product, which takes one of the copies of the holonomy along $z$
and multiplies it into the holonomy along $y$.

\begin{Truth} \lbl{CoalgebraMap} $\beta:\hat{S}(\frakg)\to\hat{U}(\frakg)$
is a co-algebra map.
\end{Truth}

\begin{AlmostTruth} \lbl{AlgebraMap}
$\beta:\hat{S}(\frakg)\to\hat{U}(\frakg)$ is an algebra map as well.
\end{AlmostTruth}

Together,~\ref{CoalgebraMap} and~\ref{AlgebraMap} say that~\eqref{Uglevel}
is also valid in $\hat{S}(\frakg)$. Identifying $\hat{S}(\frakg)$ with the
space $F(\frakg^\star)$ of functions on $\frakg^\star$, the product
$\times_S$ of $\hat{S}(\frakg)$ becomes the diagonal map $f(x,y)\mapsto
f(x,x)$, the co-product $\Delta$ becomes the map $f(x)\mapsto f(x+y)$,
and~\eqref{Uglevel} becomes
\[ h_A(L_2)(x,y,z) = h_A(L_1)(x,y,y+z). \]
But now it's clear why $\Aa_A$ is invariant under the second Kirby
move---on holonomies, the second Kirby move is just a change of variables,
which doesn't change the integral!

Let us now see why the ``almost'' in almost truth~\ref{AlgebraMap} is
harmless in our case. (If you wish to take that on faith, skip straight to
section~\ref{subsec:ThereAre}).

\begin{definition} \lbl{def:PureDiff}
A differential operator is said to be {\em pure} with respect to some
variable $y$ if its coefficients are independent of $y$ and
every term in it is of positive order in $\partial/\partial y$. We allow
infinite order differential operators, provided they are ``convergent'' in
a sense that will be made precise in the context in which they will be
used.
\end{definition}

If $D$ is pure with respect to $y$ then the fundamental theorem of
calculus shows that the integral of $Df$ with respect to $y$ vanishes
(at least when $f$ is appropriately decreasing; the $f$'s we will use
below decrease like Gaussians, which is more than enough).  With that
in mind, the following claim and exercise explain why the ``almost''
in almost truth~\ref{AlgebraMap} is harmless:

\begin{claim} \lbl{DifferByPureDiff} When the native product $\times_S$
of $\hat{S}(\frakg)$ and the product $\times_U$ it inherits from
$\hat{U}(\frakg)$ via $\beta$ are considered as products on
functions on $\frakg^\star$, they differ by some differential
operator $D'_{\BCH}: F(\frakg^\star)\otimes F(\frakg^\star) =
F(\frakg^\star\oplus\frakg^\star)\to F(\frakg^\star\oplus\frakg^\star)$,
composed with the diagonal map $\times_S$:
\[ \times_U-\times_S = \times_S\circ D'_{\BCH}. \]
\end{claim}

\begin{proof} (sketch) A hint that claim~\ref{DifferByPureDiff}
should be true is already in the more common
formulation of the Poincare-Birkhoff-Witt theorem (see
e.g.~\cite[paragraph~2.3.6]{Dixmier:EnvelopingAlgebras}), saying that
the symmetric algebra $S(\frakg)$ is isomorphic {\em as an algebra}
to the associated graded $\gr U(\frakg)$ of the universal enveloping
algebra $U(\frakg)$. This means that if $f$ and $g$ are homogeneous
polynomials, then $f\times_S g$ and $f\times_U g$ differ by lower degree
terms, and if the force is with us these lower degree terms come from
applying differential operators.  As a typical example, let's look at
the term of degree one less. If $f=\prod_\alpha a_\alpha\in S(\frakg)$
and $g=\prod_\beta b_\beta\in S(\frakg)$, then
\[ f\times_U g = f\times_S g+\frac12\sum_{\gamma,\delta}[a_\gamma,b_\delta]
  \prod_{\alpha\neq\gamma}a_\alpha\,\prod_{\beta\neq\delta}b_\beta
  =f\times_S g+\times_S\circ D_1(f,g)\quad\text{(mod lower degrees)},
\]
where (from this formula) the operator $D_1$ can be
written in terms of a basis $\{l_i\}$ of $\frakg$ as
$
  D_1=\frac12\sum_{i,j}\left(([l_i,l_j]\otimes 1\right))
  \cdot
  \left(
    \frac{\partial}{\partial l_i}\otimes\frac{\partial}{\partial l_j}
  \right)
$.
Written in full, $D'_{\BCH}$ is an infinite sum of operators $D_k$,
where each $D_k$ is of bounded degree and involves $k$ applications of
the Lie bracket $[\cdot,\cdot]$. Below it will be justified to think of
the bracket as ``small'', and hence the sum $D'_{\BCH}$ is convergent.
\end{proof}

\begin{exercise} Show that $D=(D'_{\BCH}\otimes 1)\circ(1\otimes\Delta)$
is a pure differential operator on the space
$F(\frakg^\star)\otimes F(\frakg^\star)=F(\frakg^\star\oplus\frakg^\star)$,
with respect to all variables in the second copy of $\frakg^\star$.
\end{exercise}

\begin{exercise} \lbl{ex:BCH}
(Compare with~\cite[corollary~\ref{II-cor:BCH}]{Aarhus:II})
Use the Baker-Campbell-Hausdorff formula (see
e.g.~\cite[section~V.5]{Jacobson:LieAlgebras}) to find explicit formulas
for the operators $D'_{\BCH}$ and $D$.
\end{exercise}

\ifenroute{Thought your flight was nearly over but forgot about time zone
differences? Bummer. We have some time for some more facts.}{Here are
some more relevant facts:}

\begin{fact} The \Arhus{} integral $\Aa_A$ is invariant under reversal of
orientation of any component of the link (better be that way, for the Kirby
calculus is about unoriented links).
\end{fact}

\begin{proof} On the level of $F(\frakg^\star)$, reversal of the
orientation of a component acts by negating the corresponding variable,
say, $h_A(L)(x,y,z)\to h_A(L)(x,-y,z)$. This operation does not change the
integral of $h_A(L)$.
\end{proof}

\begin{fact} $\Aa$ is independent of the choice of the base points on the
components.
\end{fact}

\begin{proof} Moving the base point on one of the components amounts to
acting on the corresponding $\hat{U}(\frakg)$ by some group element
$g$. This translates to acting on the corresponding variable of
$h_A(L)(x,y,\ldots)$ by $\Ad g$. But the adjoint action $\Ad g$ acts by
a volume preserving transformation, and hence the integral of $h_A(L)$
does not change.
\end{proof}

We leave it as an exercise to the reader to check that the normalization
$\calN$ in definition~\ref{ArhusDefinition} can be chosen so as to make $\Aa_A$
invariant under the easier ``first'' Kirby move, shown in
figure~\ref{KirbyI}. The solution appears in
section~\ref{subsubsec:Normalization}.

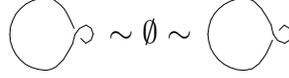
\begin{figure}[htpb]
\[ \printname{KirbyI}
	\setlength{\unitlength}{0.75\standardunitlength}
	\begin{array}{c}  \hspace{-1.7mm}
        	\raisebox{-8pt}{\begingroup\makeatletter\ifx\SetFigFont\undefined
\def\x#1#2#3#4#5#6#7\relax{\def\x{#1#2#3#4#5#6}}%
\expandafter\x\fmtname xxxxxx\relax \def\y{splain}%
\ifx\x\y   
\gdef\SetFigFont#1#2#3{%
  \ifnum #1<17\tiny\else \ifnum #1<20\small\else
  \ifnum #1<24\normalsize\else \ifnum #1<29\large\else
  \ifnum #1<34\Large\else \ifnum #1<41\LARGE\else
     \huge\fi\fi\fi\fi\fi\fi
  \csname #3\endcsname}%
\else
\gdef\SetFigFont#1#2#3{\begingroup
  \count@#1\relax \ifnum 25<\count@\count@25\fi
  \def\x{\endgroup\@setsize\SetFigFont{#2pt}}%
  \expandafter\x
    \csname \romannumeral\the\count@ pt\expandafter\endcsname
    \csname @\romannumeral\the\count@ pt\endcsname
  \csname #3\endcsname}%
\fi
\fi\endgroup
{\renewcommand{\dashlinestretch}{30}
\begin{picture}(2360,633)(0,-10)
\path(544,358)	(515.499,422.335)
	(488.096,475.782)
	(461.098,519.076)
	(433.815,552.954)
	(375.627,595.399)
	(308.000,609.000)

\path(308,609)	(248.477,603.595)
	(192.604,586.383)
	(141.702,558.683)
	(97.092,521.812)
	(60.096,477.091)
	(32.035,425.836)
	(14.229,369.366)
	(8.000,309.000)

\path(8,309)	(14.123,248.736)
	(31.672,192.516)
	(59.416,141.568)
	(96.125,97.125)
	(140.568,60.416)
	(191.516,32.672)
	(247.736,15.123)
	(308.000,9.000)

\path(308,9)	(351.601,15.859)
	(389.890,34.526)
	(452.311,95.820)
	(498.827,169.954)
	(533.000,234.000)

\path(533,234)	(542.039,271.489)
	(552.000,309.000)

\path(552,309)	(587.644,358.234)
	(645.000,384.000)

\path(645,384)	(687.090,357.019)
	(702.000,309.000)

\path(702,309)	(685.121,259.920)
	(642.000,234.000)

\path(642,234)	(600.000,275.000)

\path(2194,260)	(2165.499,195.665)
	(2138.096,142.218)
	(2111.098,98.924)
	(2083.815,65.046)
	(2025.627,22.601)
	(1958.000,9.000)

\path(1958,9)	(1898.477,14.405)
	(1842.604,31.617)
	(1791.702,59.317)
	(1747.092,96.188)
	(1710.096,140.909)
	(1682.035,192.164)
	(1664.229,248.634)
	(1658.000,309.000)

\path(1658,309)	(1664.123,369.264)
	(1681.672,425.484)
	(1709.416,476.432)
	(1746.125,520.875)
	(1790.568,557.584)
	(1841.516,585.328)
	(1897.736,602.877)
	(1958.000,609.000)

\path(1958,609)	(2001.601,602.141)
	(2039.890,583.474)
	(2102.311,522.180)
	(2148.827,448.046)
	(2183.000,384.000)

\path(2183,384)	(2192.039,346.511)
	(2202.000,309.000)

\path(2202,309)	(2237.644,259.766)
	(2295.000,234.000)

\path(2295,234)	(2337.090,260.981)
	(2352.000,309.000)

\path(2352,309)	(2335.121,358.080)
	(2292.000,384.000)

\path(2292,384)	(2250.000,343.000)

\put(833,234){\makebox(0,0)[lb]{\smash{{\mathmode{\sim\emptyset\sim}}}}}
\end{picture}
} }
        	\hspace{-1.9mm}
	\end{array}
 \]
\caption{The first Kirby move: an unknotted unlinked component of framing
$\pm 1$ can be removed.}
\lbl{KirbyI}
\end{figure}

\subsection{There aren't, but we can do without.} \lbl{subsec:ThereAre}
Rather than using a single highly non-trivial flat connection $A$,
we use the average holonomy of many not-necessarily-flat connections, with
respect to a non-existent measure, and show that the resulting averaged
holonomies have all the right properties. Namely, we replace $h_A(L)$ by
$h_{\frakg,k}(L):=\int h_B(L)d\mu_k(B)$, where $d\mu_k(B)$ is the
famed Chern-Simons~\cite{Witten:JonesPolynomial} measure on the space
of $\frakg$-connections, depending on some integer parameter $k$:
\[ 
  d\mu_k(B)=\exp\left(
    \frac{ik}{4\pi}\int_{S^3}\tr B\wedge dB+\frac23 B\wedge B\wedge B
  \right)\calD B
\]
($\calD B$ denotes the path integral measure over the space of
$\frakg$-connections). In other words, we set
\[ \Aa_\frakg(L)=\calN\int_{\frakg^\star\oplus\ldots\oplus\frakg^\star}
  \left(\int h_B(L)(x_1,\ldots,x_n)d\mu_k(B)\right)
  \,dx_1\ldots dx_n.
\]

Let us see why $h_{\frakg,k}(L)$ has all the right properties:
\begin{itemize}
\item {\em Flatness: } We only need to know that $h_{\frakg,k}(L)$
  is a link invariant. Indeed it is, by the usual topological
  invariance of the Chern-Simons path
  integral~\cite{Witten:JonesPolynomial}.
\item {\em Non-triviality: } If $h_{\frakg,k}(L)$ was a trivial link
  invariant, so would have $\Aa_\frakg$. Fortunately the Chern-Simons
  path integral isn't trivial. \ifenroute{Look up, you were offered a
  coffee. Take it and continue reading.}{}
\end{itemize}
We divide the second Kirby move into two steps: first a component is
doubled, and then we form the connected sum of one of the copies with some
third component. So two more properties are needed:
\begin{itemize}
\item {\em Behavior under doubling: } We need to know that
  if $L_2$ is obtained from $L_1$ by (say) doubling the first
  component, then $h_{\frakg,k}(L_2)=(\Delta\otimes 1\otimes\cdots\otimes
  1)h_{\frakg,k}(L_1)$. This property holds for every individual $B$, and
  it is a linear property that survives averaging.
\item {\em Behavior under connected sum: } For notational simplicity, let
  us restrict to knots. If $C_1$ and $C_2$ are two
  knots and $C_1\#C_2$ is their connected sum, we need to know that
  $h_{\frakg,k}(C_1\#C_2)=h_{\frakg,k}(C_1)\times_U h_{\frakg,k}(C_2)$.
  With the proper normalization, this is a well known property of the
  Chern-Simons path integral, and a typical example of cut-and-paste
  properties of topological quantum field theories.
\end{itemize}
We should note that in part II of this series (\cite{Aarhus:II}) we
will recombine these last two properties again into one, whose proof,
due to Le, H.~Murakami, J.~Murakami, and
Ohtsuki~\cite{LeMurakami2Ohtsuki:3Manifold}, is rather intricate and
ingenious.

\subsection{The diagrammatic case: formal Gaussian integration}
\lbl{subsec:DiagrammaticCase}
It may or may not be possible to make sense of the ideas outlined in
sections~\ref{subsec:WhatIf} and~\ref{subsec:ThereAre} as stated, per
Lie-algebra $\frakg$ and per integer $k$. We don't know though we do
want to know. Anyway, our approach is different. We do everything in
the $k\to\infty$ limit, where the Chern-Simons path integral has an
asymptotic expansion in terms of Feynman diagrams. Furthermore, as is
becoming a standard practice among topological perturbativites (see
e.g.~\cite{Bar-Natan:Thesis, Kontsevich:FeynmanDiagrams,
BottTaubes:SelfLinking, Thurston:IntegralExpressions,
AltschulerFreidel:AllOrders}), we divorce the Lie-algebras from the
Feynman diagrams and work in the universal diagrammatic setting. In
this case, the Chern-Simons path integral is valued in the space
$\calA(\circlearrowleft_X)$ of diagrams as in
figure~\ref{fig:LinkDiagrams} modulo the $STU$ relations, whose precise
form is immaterial here (figure~\ref{fig:STU}, if you insist).

\begin{figure}[htpb]
\[ \printname{LinkDiagrams}
	\setlength{\unitlength}{0.5\standardunitlength}
	\begin{array}{c}  \hspace{-1.7mm}
        	\raisebox{-8pt}{\begingroup\makeatletter\ifx\SetFigFont\undefined
\def\x#1#2#3#4#5#6#7\relax{\def\x{#1#2#3#4#5#6}}%
\expandafter\x\fmtname xxxxxx\relax \def\y{splain}%
\ifx\x\y   
\gdef\SetFigFont#1#2#3{%
  \ifnum #1<17\tiny\else \ifnum #1<20\small\else
  \ifnum #1<24\normalsize\else \ifnum #1<29\large\else
  \ifnum #1<34\Large\else \ifnum #1<41\LARGE\else
     \huge\fi\fi\fi\fi\fi\fi
  \csname #3\endcsname}%
\else
\gdef\SetFigFont#1#2#3{\begingroup
  \count@#1\relax \ifnum 25<\count@\count@25\fi
  \def\x{\endgroup\@setsize\SetFigFont{#2pt}}%
  \expandafter\x
    \csname \romannumeral\the\count@ pt\expandafter\endcsname
    \csname @\romannumeral\the\count@ pt\endcsname
  \csname #3\endcsname}%
\fi
\fi\endgroup
{\renewcommand{\dashlinestretch}{30}
\begin{picture}(4925,1247)(0,-10)
\thicklines
\put(616.000,616.000){\arc{1200.000}{3.1416}{6.2832}}
\path(6.437,863.201)(16.000,616.000)(123.892,838.619)
\put(2416.000,616.000){\arc{1200.000}{3.1416}{6.2832}}
\path(1806.437,863.201)(1816.000,616.000)(1923.892,838.619)
\put(4216.000,616.000){\arc{1200.000}{3.1416}{6.2832}}
\path(3606.437,863.201)(3616.000,616.000)(3723.892,838.619)
\put(616.000,616.000){\arc{1200.000}{6.2832}{9.4248}}
\put(2416.000,616.000){\arc{1200.000}{6.2832}{9.4248}}
\put(4216.000,616.000){\arc{1200.000}{6.2832}{9.4248}}
\thinlines
\put(628.500,653.500){\arc{285.044}{2.2318}{6.5494}}
\path(460.127,634.593)(541.000,541.000)(513.685,661.640)
\path(316,1141)(616,616)(1216,616)
\path(616,616)(241,166)
\path(1141,316)(1891,316)
\path(1141,916)(1891,916)
\path(4216,1216)(4216,16)
\path(3616,616)(3016,616)
\put(1141,91){\makebox(0,0)[lb]{\smash{{\mathmode{\scriptstyle x}}}}}
\put(2941,91){\makebox(0,0)[lb]{\smash{{\mathmode{\scriptstyle y}}}}}
\put(4741,91){\makebox(0,0)[lb]{\smash{{\mathmode{\scriptstyle z}}}}}
\end{picture}
} }
        	\hspace{-1.9mm}
	\end{array}
 \]
\caption{
  Diagrams in $\calA(\circlearrowleft_{\{x,y,z\}})$ are trivalent
  graphs made of $3$ oriented circles labeled $x$, $y$, and $z$, and
  some number of additional ``internal'' edges. ``Internal vertices'',
  in which three internal edges meet, are ``oriented'' - a cyclic order
  on the edges emanating from such a vertex is specified.}
\lbl{fig:LinkDiagrams}
\end{figure}
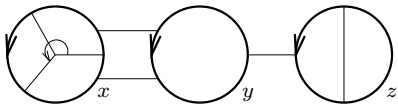

Choosing a base point on each oriented circle (akin to our choice of a
base point on each link component, in section~\ref{subsec:WhatIf}) and
cutting the circles open, we get a sum of diagrams, deserving of the
name $h_\infty(L)$, in the space $\calA(\uparrow_X)$. A typical
element of that space appears in figure~\ref{fig:PTDiagram}.

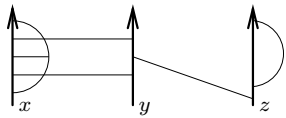
\begin{figure}[htpb]
\[ \printname{PTDiagram}
	\setlength{\unitlength}{0.5\standardunitlength}
	\begin{array}{c}  \hspace{-1.7mm}
        	\raisebox{-8pt}{\begingroup\makeatletter\ifx\SetFigFont\undefined
\def\x#1#2#3#4#5#6#7\relax{\def\x{#1#2#3#4#5#6}}%
\expandafter\x\fmtname xxxxxx\relax \def\y{splain}%
\ifx\x\y   
\gdef\SetFigFont#1#2#3{%
  \ifnum #1<17\tiny\else \ifnum #1<20\small\else
  \ifnum #1<24\normalsize\else \ifnum #1<29\large\else
  \ifnum #1<34\Large\else \ifnum #1<41\LARGE\else
     \huge\fi\fi\fi\fi\fi\fi
  \csname #3\endcsname}%
\else
\gdef\SetFigFont#1#2#3{\begingroup
  \count@#1\relax \ifnum 25<\count@\count@25\fi
  \def\x{\endgroup\@setsize\SetFigFont{#2pt}}%
  \expandafter\x
    \csname \romannumeral\the\count@ pt\expandafter\endcsname
    \csname @\romannumeral\the\count@ pt\endcsname
  \csname #3\endcsname}%
\fi
\fi\endgroup
{\renewcommand{\dashlinestretch}{30}
\begin{picture}(3407,1352)(0,-10)
\put(22.000,715.000){\arc{900.000}{4.7124}{7.8540}}
\put(2984.500,752.500){\arc{828.402}{4.8030}{7.7633}}
\thicklines
\path(22,115)(22,1315)
\path(82.000,1075.000)(22.000,1315.000)(-38.000,1075.000)
\path(1522,115)(1522,1315)
\path(1582.000,1075.000)(1522.000,1315.000)(1462.000,1075.000)
\path(3022,115)(3022,1315)
\path(3082.000,1075.000)(3022.000,1315.000)(2962.000,1075.000)
\thinlines
\path(472,715)(22,715)
\path(22,940)(1522,940)
\path(22,490)(1522,490)
\path(1522,715)(3022,190)
\put(97,40){\makebox(0,0)[lb]{\smash{{\mathmode{\scriptstyle x}}}}}
\put(1597,40){\makebox(0,0)[lb]{\smash{{\mathmode{\scriptstyle y}}}}}
\put(3097,40){\makebox(0,0)[lb]{\smash{{\mathmode{\scriptstyle z}}}}}
\end{picture}
} }
        	\hspace{-1.9mm}
	\end{array}
 \]
\caption{
  A diagram in $\calA(\uparrow_{\{x,y,z\}})$.  Here and below,
  internal trivalent vertices are always oriented counterclockwise, and
  quadrivalent vertices are just artifacts of the planar projection. This
  diagram is a cut-open of the diagram in figure~\ref{fig:LinkDiagrams}.
}
\lbl{fig:PTDiagram}
\end{figure}

\ifenroute{At this point, we are in very good shape, and things should make
sense even with the tray-tables up.}{At this point, we are in very good
shape.} There is a well known parallelism
between spaces of diagrams such as $\calA(\uparrow_X)$ and
various spaces that pertain to Lie algebras, such as $U(\frakg)^{\otimes
n}$. We have a ``holonomy'' living in the former space, and a technique
(slogan~\ref{IntegrateHolonomies}) living in the latter. All we have to do
is to imitate the Lie-level technique on the diagram level. To do that,
let us first summarize the technique of section~\ref{subsec:WhatIf} in one
line:
\[ \begin{CD}
  h_A(L)\in U(\frakg)^{\otimes n} @>PBW>> S(\frakg)^{\otimes n} @>\int >>
    \mathbb C
\end{CD}. \]

\par\noindent
{\bf The parallelism:} (The primary reference for points~\ref{item:UandAn}
to~\ref{item:CandAempty} below is~\cite{Bar-Natan:Vassiliev}. See
also\cite{Bar-Natan:Homotopy, LeMurakami:Universal}.)

\begin{enumerate}

\item \lbl{item:UandAn} The parallel of $U(\frakg)^{\otimes n}$
is, as already noted, $\calA(\uparrow_X)$. Their affinity
is first seen in the existence of a structure-respecting map
$\Tg:\calA(\uparrow_X)\to U(\frakg)^{\otimes n}$. Very briefly,
$\Tg$ is defined by mapping all internal vertices to copies of the
structure constants tensor, all internal edges the metric of $\frakg$,
and the vertical arrows to the ordered products of the Lie algebra
elements seen along them, namely to elements of $U(\frakg)$.

Notice that the parallel of the Lie bracket $[\cdot,\cdot]$ is a vertex.
Thus iterated brackets correspond to high-degree diagrams, which are
``small'' in the sense of the completed space $\calA(\uparrow_X)$.
This justifies the last sentence in the proof of
claim~\ref{DifferByPureDiff}.

\item The parallel of $S(\frakg)^{\otimes n}$ is the space $\calB(X)$
of ``$X$-marked uni-trivalent diagrams''\footnote{Called ``Chinese
characters'' in the politically incorrect~\cite{Bar-Natan:Vassiliev}.},
the space of uni-trivalent graphs whose trivalent vertices are oriented
and whose univalent vertices are marked by the symbols in the set $X$
(possibly with omissions and/or repetitions), modulo the $AS$ and $IHX$
relations, whose precise form is immaterial here
(figure~\ref{fig:IHXAS}, if you insist). An example appears in
figure~\ref{fig:TypicalCC}.  There is a structure-respecting map
$\Tg:\calB(X)\to S(\frakg)^{\otimes n}$; it maps uni-trivalent diagrams
with $k$ external legs to degree $k$ elements of $S(\frakg)^{\otimes
n}$.

\begin{figure}[htpb]
\[ \printname{TypicalCC}
	\setlength{\unitlength}{0.5\standardunitlength}
	\begin{array}{c}  \hspace{-1.7mm}
        	\raisebox{-8pt}{\begingroup\makeatletter\ifx\SetFigFont\undefined
\def\x#1#2#3#4#5#6#7\relax{\def\x{#1#2#3#4#5#6}}%
\expandafter\x\fmtname xxxxxx\relax \def\y{splain}%
\ifx\x\y   
\gdef\SetFigFont#1#2#3{%
  \ifnum #1<17\tiny\else \ifnum #1<20\small\else
  \ifnum #1<24\normalsize\else \ifnum #1<29\large\else
  \ifnum #1<34\Large\else \ifnum #1<41\LARGE\else
     \huge\fi\fi\fi\fi\fi\fi
  \csname #3\endcsname}%
\else
\gdef\SetFigFont#1#2#3{\begingroup
  \count@#1\relax \ifnum 25<\count@\count@25\fi
  \def\x{\endgroup\@setsize\SetFigFont{#2pt}}%
  \expandafter\x
    \csname \romannumeral\the\count@ pt\expandafter\endcsname
    \csname @\romannumeral\the\count@ pt\endcsname
  \csname #3\endcsname}%
\fi
\fi\endgroup
{\renewcommand{\dashlinestretch}{30}
\begin{picture}(2359,1390)(0,-10)
\put(2100,715){\ellipse{450}{450}}
\path(600,1015)(1200,1015)(1200,415)
	(600,415)(600,1015)
\path(300,1315)(600,1015)
\path(1200,1015)(1500,1315)
\path(1200,415)(1500,115)
\path(600,415)(300,115)
\path(1875,715)(2325,715)
\path(2100,940)(2100,1315)
\path(2100,490)(2100,115)
\put(0,1240){\makebox(0,0)[lb]{\smash{{\mathmode{\scriptstyle x}}}}}
\put(0,40){\makebox(0,0)[lb]{\smash{{\mathmode{\scriptstyle y}}}}}
\put(1575,40){\makebox(0,0)[lb]{\smash{{\mathmode{\scriptstyle y}}}}}
\put(1575,1240){\makebox(0,0)[lb]{\smash{{\mathmode{\scriptstyle w}}}}}
\put(2175,1240){\makebox(0,0)[lb]{\smash{{\mathmode{\scriptstyle z}}}}}
\put(2175,40){\makebox(0,0)[lb]{\smash{{\mathmode{\scriptstyle y}}}}}
\end{picture}
} }
        	\hspace{-1.9mm}
	\end{array}
 \]
\caption{An $\{x,y,z,w\}$-marked uni-trivalent diagram.}
\lbl{fig:TypicalCC}
\end{figure}
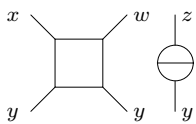

\item The parallel of the Poincare-Birkhoff-Witt isomorphism
$U(\frakg)^{\otimes n}\to S(\frakg)^{\otimes n}$ is a map
$\sigma:\calA(\uparrow_X)\to\calB(X)$, which is more easily described
through its inverse $\chi$. If $C\in \calB(X)$ is an $X$-marked
uni-trivalent diagram with $k_x$ legs marked $x$ for $x\in X$, then
$\chi(C)\in\calA(\uparrow_X)$ is the average of the $\prod_x k_x!$ ways
of attaching the legs of $C$ to $n$ labeled vertical arrows (labeled by
the elements of $X$), attaching legs marked by $x$ only to the $x$-labeled
arrow, for all $x\in X$.

\item \lbl{item:CandAempty} The parallel of $\mathbb C$ is the the
set of uni-trivalent diagrams that $\Tg$ maps to degree $0$ elements
of $S(\frakg)^{\otimes n}$ --- namely, it is the space called
$\calA(\emptyset)$ of ``manifold diagrams'' --- uni-trivalent diagrams with
no external legs. An example appears in figure~\ref{fig:VacuumCC}.

\begin{figure}[htpb]
\[ \printname{VacuumCC}
	\setlength{\unitlength}{0.5\standardunitlength}
	\begin{array}{c}  \hspace{-1.7mm}
        	\raisebox{-8pt}{\begingroup\makeatletter\ifx\SetFigFont\undefined
\def\x#1#2#3#4#5#6#7\relax{\def\x{#1#2#3#4#5#6}}%
\expandafter\x\fmtname xxxxxx\relax \def\y{splain}%
\ifx\x\y   
\gdef\SetFigFont#1#2#3{%
  \ifnum #1<17\tiny\else \ifnum #1<20\small\else
  \ifnum #1<24\normalsize\else \ifnum #1<29\large\else
  \ifnum #1<34\Large\else \ifnum #1<41\LARGE\else
     \huge\fi\fi\fi\fi\fi\fi
  \csname #3\endcsname}%
\else
\gdef\SetFigFont#1#2#3{\begingroup
  \count@#1\relax \ifnum 25<\count@\count@25\fi
  \def\x{\endgroup\@setsize\SetFigFont{#2pt}}%
  \expandafter\x
    \csname \romannumeral\the\count@ pt\expandafter\endcsname
    \csname @\romannumeral\the\count@ pt\endcsname
  \csname #3\endcsname}%
\fi
\fi\endgroup
{\renewcommand{\dashlinestretch}{30}
\begin{picture}(2420,1256)(0,-10)
\put(608,612){\ellipse{1200}{1200}}
\put(2108,612){\ellipse{600}{900}}
\path(158,987)(608,612)(1058,987)
\path(608,612)(608,12)
\path(1808,612)(2408,612)
\path(1433,987)(1883,912)
\path(1883,312)(1133,312)
\path(908,1137)	(961.376,1167.463)
	(1008.459,1191.462)
	(1050.126,1209.215)
	(1087.258,1220.944)
	(1151.428,1227.204)
	(1208.000,1212.000)

\path(1208,1212)	(1265.377,1178.903)
	(1317.052,1135.747)
	(1361.696,1084.256)
	(1397.982,1026.150)
	(1424.584,963.150)
	(1440.172,896.977)
	(1443.420,829.354)
	(1433.000,762.000)

\path(1433,762)	(1410.875,714.085)
	(1369.516,675.457)
	(1303.649,642.601)
	(1259.877,627.238)
	(1208.000,612.000)

\end{picture}
} }
        	\hspace{-1.9mm}
	\end{array}
 \]
\caption{A connected manifold diagram. All vertices are oriented
counterclockwise.}
\lbl{fig:VacuumCC}
\end{figure}
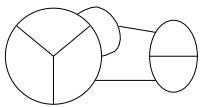

\item \lbl{item:Gint} Last, we need a parallel
$\Gint:\calB(X)\to\calA(\emptyset)$ for the partially defined integration
map $\int:S(\frakg)^{\otimes n}\to\mathbb C$. This is a new and
(hopefully) amusing ingredient, \ifenroute{and we have enough time to
elucidate on it even with the plane already rolling on the tarmac.}{so
let us say a bit more.}

\end{enumerate}

The new $\Gint$ better be defined on $\sigma h_\infty(L)$. To ensure this,
we first need some knowledge about the structure of $\sigma h_\infty(L)$,
and the following easy claim suffices:

\begin{claim} \lbl{claim:ZisPG}
If $\left(l_{xy}\right)$ is the $n\times n$ linking matrix of $L$
(that is, $l_{xy}$ is the linking number of the component $x$ with the
component $y$, and $l_{xx}$ is the self-linking of the component $x$
--- the linking number of that component with its framing), then
\begin{equation} \lbl{ImGaussian}
  \underbrace{\sigma h_\infty(L)=\exp_\udot\left(
    \frac12\sum_{x,y}l_{xy}\chordn{x}{y}+
    \left(\parbox{2.1in}{
      connected uni-trivalent diagrams that have trivalent vertices
    }\right)
  \right)}_{\text{\small\it I'm Gaussian!}}.
\end{equation}
Here $\exp_\udot$ denotes power-series exponentiation using the disjoint
union product on $\calB(X)$, and $\chordn{x}{y}$ denotes the only
connected uni-trivalent diagrams that have no trivalent vertices --- solitary
edges whose univalent ends are marked $i$ and $j$.
\end{claim}

\begin{proof} (sketch) There is a co-product
$\square:\calB(X)\to\calB(X)\otimes\calB(X)$, mapping every uni-trivalent
diagram to the sum of all possible ways of splitting it by its connected
components. A simple argument shows that $Z=\sigma h_\infty(L)$ satisfies
$\square Z=Z\otimes Z$, and thus $Z=\exp(\square\text{-primitives})$. The
$\square$-primitives of $\calB(X)$ are the connected uni-trivalent diagrams,
and hence $Z=\exp($connected characters$)$. (This is a variant of a
standard argument from quantum field theory, saying that the logarithm
of the partition function can be computed using connected Feynman
diagrams). All that is left is to determine the coefficients of the
simplest possible connected uni-trivalent diagrams, the $\chordn{x}{y}$'s.
These correspond to degree $1$ Vassiliev invariants. Namely, to linking
numbers.
\end{proof}

Equation~\eqref{ImGaussian} says ``I'm Gaussian!''. Equations don't lie, so
let's take it seriously. Remembering that $\Tg$ maps leg-count to degree
and thinking of trivalent vertices as ``small'' and hence of the second
term in~\eqref{ImGaussian} as a perturbation, we see that $\Tg\sigma
h_\infty(L)$ is indeed a Gaussian! Thus the standard Gaussian integration
technique of Feynman diagrams applies (a refresher is in
section~\ref{GaussianDummies}). But this is a diagrammatic technique, and
hence it can be applied straight at the diagrammatic level of $\sigma
h_\infty(L)$. It takes two steps:
\begin{enumerate}
\item Splitting the quadratic part out, negating and inverting it, getting
\[ \exp_\udot
  \left(
    -\frac12\sum_{x,y}l^{xy}\chordu{\partial_x}{\partial_y}
  \right),
\]
where $\left(l^{xy}\right)$ is the inverse matrix of $\left(l_{xy}\right)$
and we've introduced a new set of ``dual'' labels
$\partial_X=\{\partial_x:x\in X\}$.
\item Putting the inverted quadratic part back in and gluing its legs to
all other legs in all possible ways, making sure that the markings match.
\end{enumerate}
These steps together with all previous steps are summarized in the
commutative diagram in figure~\ref{Parallelism} and in a more graphical
form in figure~\ref{GraphicallyOriented}.

One last comment is in order. While the Chern-Simons $h_\infty$ is
perfectly good for philosophy majors, math majors find it a bit to
difficult to work with. So below they replace it by a substitute whose
properties they better understand, a variant $\Zcheck$ due
to~\cite{LeMurakami2Ohtsuki:3Manifold} of the Kontsevich
integral~\cite{Kontsevich:Vassiliev, Bar-Natan:Vassiliev} (which by
itself is a holonomy, of the Knizhnik-Zamolodchikov connection). It is
conjectured that the original Kontsevich integral is equal to
$h_\infty$.

\begin{figure}[htpb]
\[
  \begin{CD}
    \makebox[0pt][r]{\smash{$\displaystyle \begin{array}[t]{c}
       h_\infty(L) \\
       {\scriptstyle\left(\text{or }\Zcheck(L)\right)}
     \end{array}\hspace{-9pt}\in\,$}}\calA(\uparrow_X)
   @>\sigma>>
    \calB(X)
    @>\Gint>\text{(partially defined)}>
    \calA(\emptyset)\makebox[0pt][l]{$\displaystyle \,\ni\Aa(L)$}
  \\
    @V\Tg VV
    @V\Tg VV
    @V\Tg VV
  \\
    U(\frakg)^{\otimes n}
    @>PBW>>
    S(\frakg)^{\otimes n}
    @>\int >\text{(partially defined)}>
    \mathbb C
  \end{CD}
\]
\caption{
  The bottom row is where the nonsensical section~\ref{subsec:WhatIf}
  lives.  It is also were section~\ref{subsec:ThereAre} lives, both
  in the finite $k$ case (which we do not consider below) and in the
  $k\to\infty$ limit.  The top row leads to a well defined and interesting
  invariant, and is the main focus of the rest of this series of papers.
}
\lbl{Parallelism}
\end{figure}

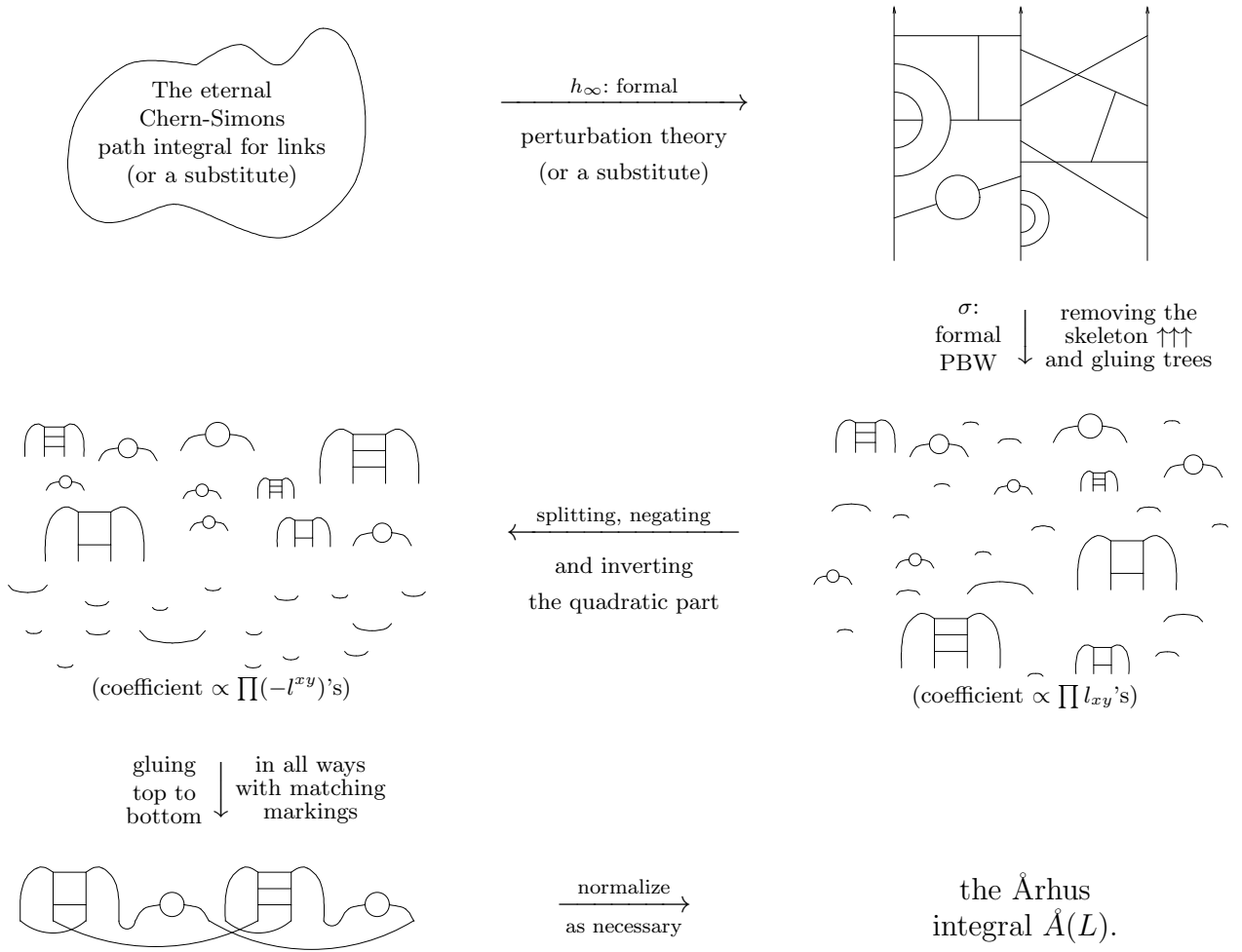
\begin{figure}[htpb]
\newcommand\TheMythical{\parbox{1.4in}{\scriptsize\begin{center}
  The eternal \\
  Chern-Simons \\
  path integral for links \\
  (or a substitute)
\end{center}}}
\[
  \begin{array}{ccc}
    \printname{MythicalCS}
	\setlength{\unitlength}{0.35\standardunitlength}
	\begin{array}{c}  \hspace{-1.7mm}
        	\raisebox{-8pt}{\begingroup\makeatletter\ifx\SetFigFont\undefined
\def\x#1#2#3#4#5#6#7\relax{\def\x{#1#2#3#4#5#6}}%
\expandafter\x\fmtname xxxxxx\relax \def\y{splain}%
\ifx\x\y   
\gdef\SetFigFont#1#2#3{%
  \ifnum #1<17\tiny\else \ifnum #1<20\small\else
  \ifnum #1<24\normalsize\else \ifnum #1<29\large\else
  \ifnum #1<34\Large\else \ifnum #1<41\LARGE\else
     \huge\fi\fi\fi\fi\fi\fi
  \csname #3\endcsname}%
\else
\gdef\SetFigFont#1#2#3{\begingroup
  \count@#1\relax \ifnum 25<\count@\count@25\fi
  \def\x{\endgroup\@setsize\SetFigFont{#2pt}}%
  \expandafter\x
    \csname \romannumeral\the\count@ pt\expandafter\endcsname
    \csname @\romannumeral\the\count@ pt\endcsname
  \csname #3\endcsname}%
\fi
\fi\endgroup
{\renewcommand{\dashlinestretch}{30}
\begin{picture}(5501,3921)(0,-10)
\path(2375,3219)	(2309.133,3229.011)
	(2245.528,3238.376)
	(2184.127,3247.096)
	(2124.873,3255.169)
	(2067.705,3262.597)
	(2012.567,3269.379)
	(1959.400,3275.514)
	(1908.145,3281.004)
	(1858.744,3285.848)
	(1811.139,3290.047)
	(1765.271,3293.599)
	(1721.082,3296.505)
	(1678.514,3298.766)
	(1637.508,3300.381)
	(1598.006,3301.350)
	(1559.949,3301.672)
	(1487.938,3300.381)
	(1421.009,3296.505)
	(1358.695,3290.047)
	(1300.528,3281.004)
	(1246.042,3269.379)
	(1194.770,3255.169)
	(1146.245,3238.376)
	(1100.000,3219.000)

\path(1100,3219)	(1028.337,3182.533)
	(954.752,3137.991)
	(879.963,3086.138)
	(842.341,3057.708)
	(804.688,3027.736)
	(767.094,2996.318)
	(729.647,2963.549)
	(655.556,2894.341)
	(583.136,2820.874)
	(547.776,2782.782)
	(513.102,2743.912)
	(479.206,2704.360)
	(446.176,2664.219)
	(414.101,2623.587)
	(383.073,2582.558)
	(353.181,2541.227)
	(324.514,2499.691)
	(297.162,2458.044)
	(271.215,2416.383)
	(246.764,2374.802)
	(223.896,2333.396)
	(202.704,2292.262)
	(183.275,2251.495)
	(165.701,2211.189)
	(150.070,2171.442)
	(136.474,2132.347)
	(125.000,2094.000)

\path(125,2094)	(113.411,2050.218)
	(102.107,2004.451)
	(91.195,1956.847)
	(80.780,1907.556)
	(70.971,1856.724)
	(61.873,1804.501)
	(53.594,1751.035)
	(46.241,1696.474)
	(39.919,1640.966)
	(34.736,1584.659)
	(30.799,1527.702)
	(28.215,1470.244)
	(27.090,1412.431)
	(27.530,1354.414)
	(29.644,1296.339)
	(33.537,1238.355)
	(39.317,1180.611)
	(47.090,1123.254)
	(56.963,1066.433)
	(69.043,1010.297)
	(83.436,954.993)
	(100.250,900.670)
	(119.591,847.476)
	(141.566,795.559)
	(166.281,745.068)
	(193.844,696.150)
	(224.361,648.955)
	(257.940,603.630)
	(294.686,560.324)
	(334.707,519.185)
	(378.109,480.361)
	(425.000,444.000)

\path(425,444)	(471.767,413.493)
	(519.659,388.335)
	(568.579,368.233)
	(618.430,352.895)
	(669.116,342.029)
	(720.540,335.340)
	(772.605,332.537)
	(825.215,333.327)
	(878.272,337.416)
	(931.682,344.514)
	(985.346,354.325)
	(1039.168,366.559)
	(1093.052,380.921)
	(1146.901,397.120)
	(1200.618,414.863)
	(1254.106,433.856)
	(1307.270,453.808)
	(1360.012,474.425)
	(1412.235,495.415)
	(1463.843,516.484)
	(1514.740,537.341)
	(1564.829,557.692)
	(1614.013,577.246)
	(1662.195,595.708)
	(1709.279,612.786)
	(1755.168,628.188)
	(1799.766,641.621)
	(1842.975,652.791)
	(1884.700,661.407)
	(1924.844,667.176)
	(1963.309,669.804)
	(2000.000,669.000)

\path(2000,669)	(2048.658,663.668)
	(2099.571,654.802)
	(2152.608,642.668)
	(2207.639,627.534)
	(2264.531,609.667)
	(2323.154,589.333)
	(2383.376,566.799)
	(2445.067,542.332)
	(2508.095,516.200)
	(2572.328,488.667)
	(2637.636,460.002)
	(2703.888,430.472)
	(2770.952,400.343)
	(2838.697,369.881)
	(2906.992,339.355)
	(2975.705,309.030)
	(3044.706,279.174)
	(3113.863,250.053)
	(3183.045,221.934)
	(3252.121,195.084)
	(3320.960,169.771)
	(3389.431,146.259)
	(3457.401,124.818)
	(3524.741,105.712)
	(3591.318,89.210)
	(3657.002,75.578)
	(3721.662,65.083)
	(3785.166,57.992)
	(3847.383,54.571)
	(3908.181,55.087)
	(3967.431,59.808)
	(4025.000,69.000)

\path(4025,69)	(4069.722,79.228)
	(4115.159,91.836)
	(4161.208,106.734)
	(4207.765,123.827)
	(4254.729,143.024)
	(4301.995,164.232)
	(4349.461,187.359)
	(4397.024,212.313)
	(4444.582,239.001)
	(4492.030,267.331)
	(4539.267,297.211)
	(4586.189,328.547)
	(4632.693,361.248)
	(4678.676,395.222)
	(4724.036,430.375)
	(4768.670,466.616)
	(4812.474,503.852)
	(4855.346,541.991)
	(4897.182,580.941)
	(4937.880,620.608)
	(4977.337,660.901)
	(5015.450,701.728)
	(5052.115,742.995)
	(5087.231,784.611)
	(5120.693,826.483)
	(5152.400,868.519)
	(5182.247,910.626)
	(5210.133,952.712)
	(5235.955,994.685)
	(5259.608,1036.452)
	(5280.991,1077.921)
	(5300.000,1119.000)

\path(5300,1119)	(5325.631,1181.914)
	(5349.128,1247.861)
	(5370.509,1316.630)
	(5389.792,1388.011)
	(5406.994,1461.795)
	(5414.820,1499.522)
	(5422.133,1537.770)
	(5428.935,1576.515)
	(5435.228,1615.728)
	(5441.014,1655.384)
	(5446.295,1695.457)
	(5451.075,1735.920)
	(5455.354,1776.748)
	(5459.135,1817.913)
	(5462.420,1859.390)
	(5465.212,1901.152)
	(5467.513,1943.173)
	(5469.324,1985.426)
	(5470.649,2027.887)
	(5471.490,2070.527)
	(5471.848,2113.321)
	(5471.726,2156.243)
	(5471.126,2199.267)
	(5470.050,2242.365)
	(5468.501,2285.513)
	(5466.481,2328.683)
	(5463.991,2371.849)
	(5461.035,2414.985)
	(5457.614,2458.065)
	(5453.731,2501.062)
	(5449.388,2543.951)
	(5444.587,2586.705)
	(5439.331,2629.297)
	(5433.621,2671.702)
	(5427.459,2713.892)
	(5420.849,2755.843)
	(5413.792,2797.527)
	(5406.290,2838.919)
	(5398.346,2879.991)
	(5389.962,2920.718)
	(5381.140,2961.074)
	(5371.882,3001.032)
	(5362.191,3040.566)
	(5352.069,3079.650)
	(5341.518,3118.257)
	(5330.540,3156.361)
	(5319.137,3193.936)
	(5295.067,3267.393)
	(5269.326,3338.419)
	(5241.930,3406.803)
	(5212.899,3472.334)
	(5182.250,3534.803)
	(5150.000,3594.000)

\path(5150,3594)	(5116.027,3645.658)
	(5072.749,3698.552)
	(5021.636,3750.042)
	(4964.158,3797.486)
	(4901.785,3838.245)
	(4835.987,3869.677)
	(4768.236,3889.142)
	(4700.000,3894.000)

\path(4700,3894)	(4635.612,3883.876)
	(4576.010,3861.236)
	(4520.536,3827.846)
	(4468.530,3785.473)
	(4419.333,3735.881)
	(4372.286,3680.836)
	(4326.729,3622.105)
	(4282.003,3561.454)
	(4237.448,3500.647)
	(4192.406,3441.452)
	(4146.218,3385.634)
	(4098.223,3334.958)
	(4047.764,3291.191)
	(3994.179,3256.098)
	(3936.811,3231.446)
	(3875.000,3219.000)

\path(3875,3219)	(3828.941,3219.450)
	(3783.796,3228.366)
	(3739.480,3244.589)
	(3695.908,3266.964)
	(3652.995,3294.332)
	(3610.655,3325.537)
	(3568.802,3359.420)
	(3527.352,3394.826)
	(3486.220,3430.597)
	(3445.319,3465.575)
	(3404.564,3498.604)
	(3363.870,3528.525)
	(3323.152,3554.183)
	(3282.325,3574.420)
	(3241.303,3588.078)
	(3200.000,3594.000)

\path(3200,3594)	(3126.050,3591.995)
	(3087.913,3587.077)
	(3048.598,3579.367)
	(3007.804,3568.726)
	(2965.227,3555.018)
	(2920.567,3538.104)
	(2873.521,3517.849)
	(2823.787,3494.113)
	(2771.062,3466.761)
	(2715.045,3435.654)
	(2655.434,3400.656)
	(2591.926,3361.629)
	(2524.219,3318.435)
	(2452.011,3270.938)
	(2414.125,3245.533)
	(2375.000,3219.000)

\put(275,1794){\makebox(0,0)[lb]{\smash{{\mathmode{\TheMythical}}}}}
\end{picture}
} }
        	\hspace{-1.9mm}
	\end{array}
 &
    \begin{CD}
      @>{h_\infty:\text{ formal}}>{\begin{array}{c}
        \text{\scriptsize perturbation theory} \\
        \text{\scriptsize (or a substitute)}
      \end{array}}>
    \end{CD} &
    \printname{BigCC}
	\setlength{\unitlength}{0.3\standardunitlength}
	\begin{array}{c}  \hspace{-1.7mm}
        	\raisebox{-8pt}{\begingroup\makeatletter\ifx\SetFigFont\undefined
\def\x#1#2#3#4#5#6#7\relax{\def\x{#1#2#3#4#5#6}}%
\expandafter\x\fmtname xxxxxx\relax \def\y{splain}%
\ifx\x\y   
\gdef\SetFigFont#1#2#3{%
  \ifnum #1<17\tiny\else \ifnum #1<20\small\else
  \ifnum #1<24\normalsize\else \ifnum #1<29\large\else
  \ifnum #1<34\Large\else \ifnum #1<41\LARGE\else
     \huge\fi\fi\fi\fi\fi\fi
  \csname #3\endcsname}%
\else
\gdef\SetFigFont#1#2#3{\begingroup
  \count@#1\relax \ifnum 25<\count@\count@25\fi
  \def\x{\endgroup\@setsize\SetFigFont{#2pt}}%
  \expandafter\x
    \csname \romannumeral\the\count@ pt\expandafter\endcsname
    \csname @\romannumeral\the\count@ pt\endcsname
  \csname #3\endcsname}%
\fi
\fi\endgroup
{\renewcommand{\dashlinestretch}{30}
\begin{picture}(5424,5439)(0,-10)
\put(12.000,3012.000){\arc{2400.000}{4.7124}{7.8540}}
\put(2712.000,912.000){\arc{600.000}{4.7124}{7.8540}}
\put(2712.000,912.000){\arc{1200.000}{4.7124}{7.8540}}
\put(12.000,3012.000){\arc{1200.000}{4.7124}{7.8540}}
\put(1362,1362){\ellipse{948}{948}}
\path(12,12)(12,5412)
\path(42.000,5292.000)(12.000,5412.000)(-18.000,5292.000)
\path(2712,12)(2712,5412)
\path(2742.000,5292.000)(2712.000,5412.000)(2682.000,5292.000)
\path(5412,12)(5412,5412)
\path(5442.000,5292.000)(5412.000,5412.000)(5382.000,5292.000)
\path(1212,3012)(2712,3012)
\path(12,4812)(2712,4812)
\path(1812,4812)(1812,3012)
\path(2712,4512)(5412,3312)
\path(5412,4812)(2712,3312)
\path(2712,2112)(5412,2112)
\path(2712,2562)(5412,912)
\path(4212,2112)(4737,3612)
\path(12,3012)(612,3012)
\path(912,1212)(12,912)
\path(1812,1512)(2712,1812)
\end{picture}
} }
        	\hspace{-1.9mm}
	\end{array}

  \\ & & \\
    & &
    \begin{CD}
      @V{
        \parbox{0.5in}{\scriptsize\begin{center}
          $\sigma$: \\ formal \\ PBW
        \end{center}}
      }V{
        \parbox{1.0in}{\scriptsize\begin{center}
          removing the \\
          skeleton $\uparrow\uparrow\uparrow$ \\
          and gluing trees
        \end{center}}
      }V
    \end{CD}
  \\ & & \\
    \begin{array}{c}
      \printname{SplitCloud}
	\setlength{\unitlength}{0.35\standardunitlength}
	\begin{array}{c}  \hspace{-1.7mm}
        	\raisebox{-8pt}{\input draws/SplitCloud.tex }
        	\hspace{-1.9mm}
	\end{array}
 \\
      \text{\scriptsize(coefficient $\propto\prod (-l^{xy})$'s)}
    \end{array} &
    \begin{CD}
      @<{\text{splitting, negating}}<{\begin{array}{c}
        \text{\scriptsize and inverting} \\
        \text{\scriptsize the quadratic part}
      \end{array}}<
    \end{CD} &
    \begin{array}{c}
      \printname{Cloud}
	\setlength{\unitlength}{0.35\standardunitlength}
	\begin{array}{c}  \hspace{-1.7mm}
        	\raisebox{-8pt}{\input draws/Cloud.tex }
        	\hspace{-1.9mm}
	\end{array}
 \\
      \text{\scriptsize(coefficient $\propto\prod l_{xy}$'s)}
    \end{array}
  \\ & & \\
    \begin{CD}
      @V{
        \parbox{0.5in}{\scriptsize\begin{center}
          gluing \\ top to \\ bottom
        \end{center}}
      }V{
        \parbox{0.85in}{\scriptsize\begin{center}
          in all ways with matching markings
        \end{center}}
      }V
    \end{CD}
  \\ & & \\
    \printname{GluedCloud}
	\setlength{\unitlength}{0.35\standardunitlength}
	\begin{array}{c}  \hspace{-1.7mm}
        	\raisebox{-8pt}{\begingroup\makeatletter\ifx\SetFigFont\undefined
\def\x#1#2#3#4#5#6#7\relax{\def\x{#1#2#3#4#5#6}}%
\expandafter\x\fmtname xxxxxx\relax \def\y{splain}%
\ifx\x\y   
\gdef\SetFigFont#1#2#3{%
  \ifnum #1<17\tiny\else \ifnum #1<20\small\else
  \ifnum #1<24\normalsize\else \ifnum #1<29\large\else
  \ifnum #1<34\Large\else \ifnum #1<41\LARGE\else
     \huge\fi\fi\fi\fi\fi\fi
  \csname #3\endcsname}%
\else
\gdef\SetFigFont#1#2#3{\begingroup
  \count@#1\relax \ifnum 25<\count@\count@25\fi
  \def\x{\endgroup\@setsize\SetFigFont{#2pt}}%
  \expandafter\x
    \csname \romannumeral\the\count@ pt\expandafter\endcsname
    \csname @\romannumeral\the\count@ pt\endcsname
  \csname #3\endcsname}%
\fi
\fi\endgroup
{\renewcommand{\dashlinestretch}{30}
\begin{picture}(7230,1562)(0,-10)
\path(622,1436)(1222,1436)(1222,836)
	(622,836)(622,1436)
\path(622,836)(622,536)
\path(1222,836)(1222,536)
\path(622,1436)	(569.605,1469.531)
	(523.223,1495.721)
	(481.976,1514.790)
	(444.985,1526.956)
	(380.254,1531.463)
	(322.000,1511.000)

\path(322,1511)	(276.753,1483.062)
	(235.552,1452.390)
	(164.847,1381.416)
	(135.123,1340.401)
	(109.006,1295.224)
	(86.385,1245.528)
	(67.150,1190.956)
	(51.192,1131.151)
	(38.400,1065.756)
	(28.666,994.414)
	(24.910,956.401)
	(21.878,916.768)
	(19.555,875.469)
	(17.928,832.460)
	(16.982,787.697)
	(16.704,741.135)
	(17.081,692.728)
	(18.098,642.434)
	(19.743,590.206)
	(22.000,536.000)

\path(1222,1436)	(1274.395,1469.531)
	(1320.777,1495.721)
	(1362.024,1514.790)
	(1399.015,1526.956)
	(1463.746,1531.463)
	(1522.000,1511.000)

\path(1522,1511)	(1567.249,1483.063)
	(1608.451,1452.393)
	(1679.159,1381.422)
	(1708.883,1340.408)
	(1735.001,1295.231)
	(1757.623,1245.536)
	(1776.858,1190.964)
	(1792.816,1131.159)
	(1805.607,1065.763)
	(1815.341,994.421)
	(1819.096,956.408)
	(1822.128,916.774)
	(1824.450,875.474)
	(1826.077,832.465)
	(1827.022,787.701)
	(1827.299,741.138)
	(1826.922,692.731)
	(1825.903,642.435)
	(1824.258,590.207)
	(1822.000,536.000)

\put(2797,836){\ellipse{450}{450}}
\path(2572,836)	(2510.982,828.237)
	(2458.441,820.452)
	(2413.498,812.427)
	(2375.275,803.941)
	(2315.472,784.710)
	(2272.000,761.000)

\path(2272,761)	(2206.262,688.576)
	(2168.205,624.944)
	(2146.340,583.965)
	(2122.000,536.000)

\path(3022,836)	(3083.018,828.237)
	(3135.559,820.452)
	(3180.502,812.427)
	(3218.725,803.941)
	(3278.528,784.710)
	(3322.000,761.000)

\path(3322,761)	(3387.737,688.576)
	(3425.795,624.944)
	(3447.660,583.965)
	(3472.000,536.000)

\put(6547,836){\ellipse{450}{450}}
\path(6322,836)	(6260.982,828.237)
	(6208.441,820.452)
	(6163.498,812.427)
	(6125.275,803.941)
	(6065.472,784.710)
	(6022.000,761.000)

\path(6022,761)	(5956.263,688.576)
	(5918.205,624.944)
	(5896.340,583.965)
	(5872.000,536.000)

\path(6772,836)	(6833.018,828.237)
	(6885.559,820.452)
	(6930.502,812.427)
	(6968.725,803.941)
	(7028.528,784.710)
	(7072.000,761.000)

\path(7072,761)	(7137.737,688.576)
	(7175.795,624.944)
	(7197.660,583.965)
	(7222.000,536.000)

\path(4373,1435)(4973,1435)(4973,835)
	(4373,835)(4373,1435)
\path(4373,835)(4373,535)
\path(4973,835)(4973,535)
\path(4372,1136)(4972,1136)
\path(4373,1435)	(4320.605,1468.531)
	(4274.223,1494.721)
	(4232.976,1513.790)
	(4195.985,1525.956)
	(4131.254,1530.463)
	(4073.000,1510.000)

\path(4073,1510)	(4027.753,1482.062)
	(3986.552,1451.390)
	(3915.847,1380.416)
	(3886.123,1339.401)
	(3860.006,1294.224)
	(3837.385,1244.528)
	(3818.150,1189.956)
	(3802.192,1130.151)
	(3789.400,1064.756)
	(3779.666,993.414)
	(3775.910,955.401)
	(3772.878,915.768)
	(3770.555,874.469)
	(3768.928,831.460)
	(3767.982,786.697)
	(3767.704,740.135)
	(3768.081,691.728)
	(3769.098,641.434)
	(3770.743,589.206)
	(3773.000,535.000)

\path(4973,1435)	(5025.395,1468.531)
	(5071.777,1494.721)
	(5113.024,1513.790)
	(5150.015,1525.956)
	(5214.746,1530.463)
	(5273.000,1510.000)

\path(5273,1510)	(5318.249,1482.063)
	(5359.451,1451.393)
	(5430.159,1380.422)
	(5459.883,1339.408)
	(5486.001,1294.231)
	(5508.623,1244.536)
	(5527.858,1189.964)
	(5543.816,1130.159)
	(5556.607,1064.763)
	(5566.341,993.421)
	(5570.096,955.408)
	(5573.128,915.774)
	(5575.450,874.474)
	(5577.077,831.465)
	(5578.022,786.701)
	(5578.299,740.138)
	(5577.922,691.731)
	(5576.903,641.435)
	(5575.258,589.207)
	(5573.000,535.000)

\put(621.000,1224.500){\arc{1825.000}{0.8533}{2.2883}}
\put(1971.000,537.000){\arc{300.000}{6.2832}{9.4248}}
\put(5721.000,612.000){\arc{335.410}{0.4636}{2.6779}}
\put(5346.000,3601.286){\arc{7184.838}{1.0217}{2.1199}}
\put(2496.000,3993.250){\arc{7864.169}{1.0737}{2.0679}}
\put(4371.000,1224.500){\arc{1825.000}{0.8533}{2.2883}}
\end{picture}
} }
        	\hspace{-1.9mm}
	\end{array}
 &
    \begin{CD} @>{\text{normalize}}>{\text{as necessary}}> \end{CD} &
    \parbox{2in}{\begin{center}
      the \Arhus{} \\ integral $\Aa(L)$.
    \end{center}}
  \end{array}
\]
\caption{
  Introduction for the graphically oriented. All diagrams here are
  representatives of big sums of diagrams.
} 
\lbl{GraphicallyOriented}
\end{figure}

\section{Introduction for math majors}
\label{sec:intro2}
After that flight of fancy, it's time to get down to earth. In this
introduction, our goal is very modest: to give the precise definition of
an invariant $\Aa$ of rational homology spheres, and to state its main
properties. The proofs of these properties of $\Aa$ (and even the fact
that it is well defined), though mostly natural and conceptual, will be
postponed to the later parts of this series.

\subsection{The definition of the \Arhus{} integral} All (3-dimensional)
rational homology spheres are surgeries on algebraically regular framed
links in $S^3$ (``regular links'' throughout this series, precise
definition below), and all regular links are closures of algebraically
regular framed pure tangles (``regular pure tangles'', precise definition
below). All that we do in this section is to define a certain invariant
$\Aa$ of regular pure tangles, following the philosophers' ideas
(section~\ref{sec:intro1}). In part II of this series
(\cite{Aarhus:II}) we follow the same philosophers' ideas and show that
they lead to simple proofs that $\Aa$ descends to an invariant of
regular links and then to an invariant of rational homology spheres.

\subsubsection{Domain and target spaces.} Let us start with the
definitions of the domain and target spaces of $\Aa$:

\par\noindent
\parbox{5.1in}{
  \begin{definition} An ($X$-marked) pure tangle (also called ``string
  link'') is an embedding $T$ of $n$ copies of the unit interval,
  $I\times\{1\},\ldots,I\times\{n\}$ into $I\times{\mathbb C}$, so that
  $T((\epsilon,i))=(\epsilon,i)$ for all $\epsilon\in\{0,1\}$ and $1\leq
  i\leq n$, considered up to endpoint-preserving isotopies. We
  assume that the components of the pure tangle are labeled by
  labels in some $n$-element label set $X$. An examples is on the right.
  \end{definition}
}\quad
$\printname{PureTangles}
	\setlength{\unitlength}{0.5\standardunitlength}
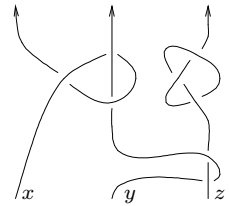
	\begin{array}{c}  \hspace{-1.7mm}
        	\raisebox{-8pt}{\begingroup\makeatletter\ifx\SetFigFont\undefined
\def\x#1#2#3#4#5#6#7\relax{\def\x{#1#2#3#4#5#6}}%
\expandafter\x\fmtname xxxxxx\relax \def\y{splain}%
\ifx\x\y   
\gdef\SetFigFont#1#2#3{%
  \ifnum #1<17\tiny\else \ifnum #1<20\small\else
  \ifnum #1<24\normalsize\else \ifnum #1<29\large\else
  \ifnum #1<34\Large\else \ifnum #1<41\LARGE\else
     \huge\fi\fi\fi\fi\fi\fi
  \csname #3\endcsname}%
\else
\gdef\SetFigFont#1#2#3{\begingroup
  \count@#1\relax \ifnum 25<\count@\count@25\fi
  \def\x{\endgroup\@setsize\SetFigFont{#2pt}}%
  \expandafter\x
    \csname \romannumeral\the\count@ pt\expandafter\endcsname
    \csname @\romannumeral\the\count@ pt\endcsname
  \csname #3\endcsname}%
\fi
\fi\endgroup
{\renewcommand{\dashlinestretch}{30}
\begin{picture}(2667,2467)(0,-10)
\path(1238.000,2320.000)(1208.000,2440.000)(1178.000,2320.000)
\path(1208,2440)(1208,1315)
\path(2408,40)(2408,490)
\path(8,40)	(19.388,81.284)
	(30.641,121.824)
	(41.763,161.629)
	(52.757,200.706)
	(63.626,239.066)
	(74.375,276.716)
	(95.521,349.922)
	(116.225,420.392)
	(136.512,488.196)
	(156.411,553.401)
	(175.949,616.077)
	(195.154,676.292)
	(214.053,734.115)
	(232.673,789.615)
	(251.042,842.860)
	(269.188,893.919)
	(287.137,942.860)
	(304.918,989.753)
	(322.558,1034.665)
	(340.083,1077.666)
	(357.523,1118.824)
	(374.904,1158.208)
	(392.253,1195.887)
	(426.968,1266.403)
	(461.887,1330.921)
	(497.230,1389.990)
	(533.217,1444.160)
	(570.067,1493.980)
	(608.000,1540.000)

\path(608,1540)	(677.945,1604.976)
	(724.398,1637.773)
	(780.556,1671.948)
	(847.958,1708.379)
	(886.355,1727.715)
	(928.140,1747.945)
	(973.505,1769.179)
	(1022.641,1791.526)
	(1075.742,1815.096)
	(1133.000,1840.000)

\path(1283,1840)	(1343.487,1795.860)
	(1393.409,1755.376)
	(1464.196,1681.855)
	(1500.635,1612.407)
	(1508.000,1540.000)

\path(1508,1540)	(1497.488,1484.148)
	(1475.413,1432.077)
	(1443.699,1384.518)
	(1404.275,1342.202)
	(1359.066,1305.862)
	(1310.000,1276.227)
	(1259.002,1254.029)
	(1208.000,1240.000)

\path(1208,1240)	(1157.094,1234.176)
	(1105.059,1236.733)
	(1050.357,1248.331)
	(991.449,1269.629)
	(926.797,1301.285)
	(854.864,1343.960)
	(815.686,1369.635)
	(774.111,1398.312)
	(729.946,1430.073)
	(683.000,1465.000)

\path(533,1615)	(482.400,1647.994)
	(435.566,1679.073)
	(392.333,1708.375)
	(352.536,1736.035)
	(282.591,1786.985)
	(224.412,1833.021)
	(176.683,1875.242)
	(138.083,1914.746)
	(107.295,1952.633)
	(83.000,1990.000)

\path(83,1990)	(54.093,2056.631)
	(42.728,2098.118)
	(33.129,2146.810)
	(25.076,2204.026)
	(18.350,2271.084)
	(15.416,2308.716)
	(12.732,2349.303)
	(10.269,2393.009)
	(8.000,2440.000)

\path(43.339,2321.462)(8.000,2440.000)(-16.601,2318.778)
\path(1208,1165)	(1206.432,1119.630)
	(1205.311,1080.505)
	(1204.415,1018.356)
	(1205.311,973.279)
	(1208.000,940.000)

\path(1208,940)	(1212.294,872.390)
	(1214.729,831.224)
	(1219.130,788.091)
	(1226.827,745.197)
	(1239.151,704.747)
	(1283.000,640.000)

\path(1283,640)	(1340.098,601.212)
	(1402.731,574.183)
	(1470.092,557.367)
	(1541.372,549.219)
	(1615.761,548.193)
	(1653.869,549.868)
	(1692.452,552.744)
	(1731.407,556.628)
	(1770.635,561.326)
	(1810.033,566.647)
	(1849.501,572.395)
	(1888.938,578.379)
	(1928.243,584.404)
	(1967.314,590.279)
	(2006.051,595.809)
	(2044.352,600.802)
	(2082.116,605.064)
	(2155.630,610.623)
	(2225.784,610.942)
	(2291.770,604.475)
	(2352.778,589.676)
	(2408.000,565.000)

\path(2408,565)	(2460.851,526.923)
	(2511.492,473.616)
	(2547.888,409.751)
	(2558.000,340.000)

\path(2558,340)	(2539.524,302.669)
	(2483.000,265.000)

\path(2333,265)	(2278.582,271.788)
	(2226.038,278.021)
	(2175.321,283.696)
	(2126.382,288.809)
	(2079.174,293.357)
	(2033.648,297.336)
	(1989.756,300.744)
	(1947.450,303.576)
	(1906.683,305.829)
	(1867.405,307.500)
	(1793.128,309.082)
	(1724.235,308.293)
	(1660.340,305.106)
	(1601.060,299.495)
	(1546.009,291.430)
	(1494.804,280.886)
	(1447.060,267.834)
	(1402.392,252.246)
	(1360.415,234.097)
	(1320.746,213.357)
	(1283.000,190.000)

\path(1283,190)	(1243.074,141.460)
	(1225.809,99.103)
	(1208.000,40.000)

\path(2408,640)	(2414.723,714.659)
	(2419.526,779.182)
	(2422.407,834.668)
	(2423.367,882.216)
	(2422.407,922.925)
	(2419.526,957.893)
	(2408.000,1015.000)

\path(2408,1015)	(2381.057,1070.399)
	(2338.426,1132.679)
	(2293.083,1192.369)
	(2258.000,1240.000)

\path(2258,1240)	(2235.966,1273.048)
	(2205.965,1318.049)
	(2164.482,1380.275)
	(2138.335,1419.495)
	(2108.000,1465.000)

\path(2183,1240)	(2116.732,1221.842)
	(2060.441,1208.871)
	(2013.248,1201.089)
	(1974.275,1198.495)
	(1917.472,1208.871)
	(1883.000,1240.000)

\path(1883,1240)	(1879.508,1312.342)
	(1897.894,1353.669)
	(1924.430,1396.037)
	(1955.162,1437.651)
	(1986.137,1476.714)
	(2033.000,1540.000)

\path(2033,1540)	(2053.723,1574.010)
	(2083.288,1619.331)
	(2125.208,1681.237)
	(2151.900,1720.057)
	(2183.000,1765.000)

\path(2333,1915)	(2364.506,1969.198)
	(2386.336,2009.920)
	(2408.000,2065.000)

\path(2408,2065)	(2411.527,2098.571)
	(2411.135,2140.465)
	(2408.000,2215.000)

\path(2408,2215)	(2408.000,2248.048)
	(2408.000,2293.049)
	(2408.000,2355.275)
	(2408.000,2394.495)
	(2408.000,2440.000)

\path(2438.000,2320.000)(2408.000,2440.000)(2378.000,2320.000)
\path(1958,1615)	(1923.183,1676.291)
	(1896.074,1729.026)
	(1876.454,1774.086)
	(1864.104,1812.348)
	(1860.332,1871.995)
	(1883.000,1915.000)

\path(1883,1915)	(1917.922,1937.009)
	(1962.419,1944.002)
	(2013.411,1939.004)
	(2067.819,1925.039)
	(2122.561,1905.131)
	(2174.559,1882.306)
	(2220.732,1859.587)
	(2258.000,1840.000)

\path(2258,1840)	(2297.890,1819.427)
	(2344.548,1795.095)
	(2394.248,1766.466)
	(2443.265,1733.001)
	(2487.874,1694.161)
	(2524.349,1649.407)
	(2548.966,1598.199)
	(2558.000,1540.000)

\path(2558,1540)	(2546.574,1476.369)
	(2508.781,1421.174)
	(2439.348,1369.142)
	(2391.118,1342.664)
	(2333.000,1315.000)

\put(83,40){\makebox(0,0)[lb]{\smash{{\mathmode{\scriptstyle x}}}}}
\put(1358,40){\makebox(0,0)[lb]{\smash{{\mathmode{\scriptstyle y}}}}}
\put(2483,40){\makebox(0,0)[lb]{\smash{{\mathmode{\scriptstyle z}}}}}
\end{picture}
} }
        	\hspace{-1.9mm}
	\end{array}
$

Similarly to many other knotted objects, pure tangles can be ``framed'' (we
omit the precise definition), and similarly to braids, framed pure tangles
can be closed to form a framed link. By pullback, this allows one to define
linking numbers and self-linking numbers for framed pure tangles.

\begin{definition} \lbl{def:RPT}
A framed link is called ``algebraically regular'' (``regular link'', in
short) if its linking matrix (with self-linkings on the diagonal) is
invertible. A framed tangle is called ``algebraically regular''
(``regular pure tangle'', in short), if the same condition holds, or,
alternatively, if its closure is a regular link. Let $RPT$ be the set
of all ($X$-marked) regular pure tangles.
\end{definition}

This completes the definition of the domain space of $\Aa$. Let's turn to
the target space:

\begin{definition} \lbl{def:ManifoldDiagrams}
A ``manifold diagram'' is a trivalent graph with oriented vertices
(i.e., a cyclic order is specified on the edges emanating from each
vertex), such as the one in figure~\ref{fig:VacuumCC}. The ``degree''
of a manifold diagram is half the number of vertices it has. The target
space of $\Aa$, called $\calA(\emptyset)$, is the graded completion of
the linear space spanned by all manifold diagrams modulo the $IHX$ and
$AS$ (antisymmetry) relations displayed in figure~\ref{fig:IHXAS}.
\end{definition}

\begin{figure}[htpb]
\[ IHX:\quad\printname{IHX}
	\setlength{\unitlength}{0.5\standardunitlength}
	\begin{array}{c}  \hspace{-1.7mm}
        	\raisebox{-8pt}{\begingroup\makeatletter\ifx\SetFigFont\undefined
\def\x#1#2#3#4#5#6#7\relax{\def\x{#1#2#3#4#5#6}}%
\expandafter\x\fmtname xxxxxx\relax \def\y{splain}%
\ifx\x\y   
\gdef\SetFigFont#1#2#3{%
  \ifnum #1<17\tiny\else \ifnum #1<20\small\else
  \ifnum #1<24\normalsize\else \ifnum #1<29\large\else
  \ifnum #1<34\Large\else \ifnum #1<41\LARGE\else
     \huge\fi\fi\fi\fi\fi\fi
  \csname #3\endcsname}%
\else
\gdef\SetFigFont#1#2#3{\begingroup
  \count@#1\relax \ifnum 25<\count@\count@25\fi
  \def\x{\endgroup\@setsize\SetFigFont{#2pt}}%
  \expandafter\x
    \csname \romannumeral\the\count@ pt\expandafter\endcsname
    \csname @\romannumeral\the\count@ pt\endcsname
  \csname #3\endcsname}%
\fi
\fi\endgroup
{\renewcommand{\dashlinestretch}{30}
\begin{picture}(2724,639)(0,-10)
\put(1738,235){\makebox(0,0)[lb]{\smash{{\mathmode{-}}}}}
\path(12,612)(612,612)
\path(312,612)(312,12)
\path(12,12)(612,12)
\path(1062,612)(1062,12)
\path(1062,312)(1662,312)
\path(1662,612)(1662,12)
\path(2112,612)(2712,12)
\path(2712,612)(2112,12)
\path(2262,162)(2562,162)
\put(687,237){\makebox(0,0)[lb]{\smash{{\mathmode{=}}}}}
\end{picture}
} }
        	\hspace{-1.9mm}
	\end{array}
\qquad\qquad AS:\quad\printname{AS}
	\setlength{\unitlength}{0.5\standardunitlength}
	\begin{array}{c}  \hspace{-1.7mm}
        	\raisebox{-8pt}{\begingroup\makeatletter\ifx\SetFigFont\undefined
\def\x#1#2#3#4#5#6#7\relax{\def\x{#1#2#3#4#5#6}}%
\expandafter\x\fmtname xxxxxx\relax \def\y{splain}%
\ifx\x\y   
\gdef\SetFigFont#1#2#3{%
  \ifnum #1<17\tiny\else \ifnum #1<20\small\else
  \ifnum #1<24\normalsize\else \ifnum #1<29\large\else
  \ifnum #1<34\Large\else \ifnum #1<41\LARGE\else
     \huge\fi\fi\fi\fi\fi\fi
  \csname #3\endcsname}%
\else
\gdef\SetFigFont#1#2#3{\begingroup
  \count@#1\relax \ifnum 25<\count@\count@25\fi
  \def\x{\endgroup\@setsize\SetFigFont{#2pt}}%
  \expandafter\x
    \csname \romannumeral\the\count@ pt\expandafter\endcsname
    \csname @\romannumeral\the\count@ pt\endcsname
  \csname #3\endcsname}%
\fi
\fi\endgroup
{\renewcommand{\dashlinestretch}{30}
\begin{picture}(2029,639)(0,-10)
\path(12,612)(312,312)(612,612)
\path(312,312)(312,12)
\path(1362,312)(1362,12)
\path(1062,612)	(1110.175,601.033)
	(1154.410,590.038)
	(1194.844,578.933)
	(1231.613,567.634)
	(1294.707,544.128)
	(1344.791,518.861)
	(1410.324,460.408)
	(1437.000,387.000)

\path(1437,387)	(1416.615,336.139)
	(1362.000,312.000)

\path(1362,312)	(1307.385,336.139)
	(1287.000,387.000)

\path(1287,387)	(1313.670,460.408)
	(1379.201,518.861)
	(1429.286,544.128)
	(1492.381,567.634)
	(1529.152,578.933)
	(1569.587,590.038)
	(1613.824,601.033)
	(1662.000,612.000)

\put(687,237){\makebox(0,0)[lb]{\smash{{\mathmode{+}}}}}
\put(1737,237){\makebox(0,0)[lb]{\smash{{\mathmode{=\ 0}}}}}
\end{picture}
} }
        	\hspace{-1.9mm}
	\end{array}
 \]
\caption{The $IHX$ and $AS$ relations. Each equality here represents a
whole family of relations, obtained by completing each of the diagram-stubs
shown to a full diagram in all possible ways, but in the same way within
each equality.}
\lbl{fig:IHXAS}
\end{figure}

\subsubsection{Intermediate spaces.} Our map $\Aa:RPT\to\calA(\emptyset)$
is a composition of several maps. Let us now define the intermediate
spaces we pass through:

\begin{definition} \lbl{def:TangleDiagrams}
An ``$X$-marked pure tangle diagram'' is a graph $D$ made of the following
types of edges and vertices:
\begin{itemize}
\item Edges: $n$ upward-pointing vertical directed lines marked by the
  elements of $X$ (whose union is
  ``the skeleton of $D$''), and some number of undirected edges, sometimes
  called ``chords'' or ``internal edges''.
\item Vertices: the endpoints of the skeleton, vertices in which an
  internal edge ends on the skeleton, and oriented trivalent vertices
  in which three internal edges meet.
\end{itemize}
The graph $D$ should be ``connected modulo its skeleton''. Namely,
if the skeleton of $D$ is collapsed to a single point, the
resulting graph should be connected. An example appears in
figure~\ref{fig:PTDiagram}. The ``degree'' of $D$ is half the number
of trivalent vertices it has. The graded completion of the space of all
$X$-marked pure tangle diagram modulo the $STU$ relations displayed in
figure~\ref{fig:STU} is denoted by $\calA(\uparrow_X)$. If $X=\{x\}$
is a singleton, we set $\calA=\calA(\uparrow_x)=\calA(\uparrow_X)$. We
note that the $STU$ relations implies the $IHX$ and $AS$ relations,
see~\cite{Bar-Natan:Vassiliev}.
\end{definition}

\begin{figure}[htpb]
\[ \printname{STU}
	\setlength{\unitlength}{0.5\standardunitlength}
	\begin{array}{c}  \hspace{-1.7mm}
        	\raisebox{-8pt}{\begingroup\makeatletter\ifx\SetFigFont\undefined
\def\x#1#2#3#4#5#6#7\relax{\def\x{#1#2#3#4#5#6}}%
\expandafter\x\fmtname xxxxxx\relax \def\y{splain}%
\ifx\x\y   
\gdef\SetFigFont#1#2#3{%
  \ifnum #1<17\tiny\else \ifnum #1<20\small\else
  \ifnum #1<24\normalsize\else \ifnum #1<29\large\else
  \ifnum #1<34\Large\else \ifnum #1<41\LARGE\else
     \huge\fi\fi\fi\fi\fi\fi
  \csname #3\endcsname}%
\else
\gdef\SetFigFont#1#2#3{\begingroup
  \count@#1\relax \ifnum 25<\count@\count@25\fi
  \def\x{\endgroup\@setsize\SetFigFont{#2pt}}%
  \expandafter\x
    \csname \romannumeral\the\count@ pt\expandafter\endcsname
    \csname @\romannumeral\the\count@ pt\endcsname
  \csname #3\endcsname}%
\fi
\fi\endgroup
{\renewcommand{\dashlinestretch}{30}
\begin{picture}(2724,639)(0,-10)
\put(1738,236){\makebox(0,0)[lb]{\smash{{\mathmode{-}}}}}
\path(612,12)(612,612)
\path(642.000,492.000)(612.000,612.000)(582.000,492.000)
\path(1662,12)(1662,612)
\path(1692.000,492.000)(1662.000,612.000)(1632.000,492.000)
\path(2712,12)(2712,612)
\path(2742.000,492.000)(2712.000,612.000)(2682.000,492.000)
\path(12,612)(387,312)(12,12)
\path(387,312)(612,312)
\path(1062,612)(1662,462)
\path(1062,12)(1662,162)
\path(2112,612)(2712,162)
\path(2112,12)(2712,462)
\put(687,237){\makebox(0,0)[lb]{\smash{{\mathmode{=}}}}}
\end{picture}
} }
        	\hspace{-1.9mm}
	\end{array}
 \]
\caption{The $STU$ family of relations with only diagram-stubs shown.}
\lbl{fig:STU}
\end{figure}

\begin{definition} \lbl{def:UnitrivalentDiagrams}
An ``$X$-marked uni-trivalent diagram'' is a graph $C$ made of undirected
edges and two types of vertices: oriented trivalent vertices
(``internal vertices'') and univalent vertices marked by elements of the
label set $X$ (the ``legs'' of $C$). The graph $C$
should be connected modulo its univalent vertices. Namely, if the
univalent vertices of $C$ are all joined, the resulting graph should be
connected. An example appears in figure~\ref{fig:TypicalCC}. The
``degree'' of $D$ is half the total number of vertices it has. The graded
completion of the space of all $X$-marked uni-trivalent diagrams modulo
the $IHX$ and $AS$ relations of figure~\ref{fig:IHXAS} is denoted by
$\calB(X)$. The space $\calB(X)$ is an algebra with the bilinear extension
of the disjoint union operation $\udot$ as a
\if\enroute y{
  pro\hspace{-2pt}\raisebox{-3pt}{duct}\footnote{
    The jag in this word is intentional, to compensate for a bump on the
    road. If you plan to re-read this paper on a sound surface, change the
    first line of the \LaTeX{} source to 
    $\backslash\text{\tt def}\backslash\text{\tt enroute}\{\text{\tt n}\}$
    and re-print.
  }.
} \else {
  product.
} \fi
\end{definition}

\subsubsection{Maps.} The pre-normalized \Arhus{} integral $\Aa_0$ is the
following composition:
\begin{equation} \lbl{eq:ArhusComposition} \begin{CD}
  \Aa_0:\left\{\parbox{0.5in}{\begin{center}
    regular pure tangles
  \end{center}}\right\}=RPT
  @>\Zcheck>{\parbox{1.1in}{\scriptsize\begin{center}
    the \cite{LeMurakami2Ohtsuki:3Manifold} version of the Kontsevich integral
  \end{center}}}>
  \calA(\uparrow_X)
  @>\sigma>{\parbox{0.4in}{\scriptsize\begin{center}
    formal PBW
  \end{center}}}>
  \calB(X)
  @>\Gint>{\parbox{0.65in}{\scriptsize\begin{center}
    formal Gaussian integration
  \end{center}}}>
  \calA(\emptyset).
\end{CD} \end{equation}

We just have to recall the definitions of the maps $\Zcheck$ and $\sigma$,
and define the map $\Gint$.

\begin{definition} \lbl{def:Zcheck}
The map $\Zcheck$ was defined by Le, H.~Murakami, J.~Murakami,
and Ohtsuki in~\cite{LeMurakami2Ohtsuki:3Manifold}. It is
the usual framed version of the Kontsevich integral $Z$
(see~\cite{LeMurakami:Universal}, or a simpler definition
in~\cite[section~2.2]{Bar-NatanGaroufalidisRozanskyThurston:WheelsWheeling}\footnote{
  Strictly
  speaking,~\cite{Bar-NatanGaroufalidisRozanskyThurston:WheelsWheeling}
  deals only with knots. But the generalization to links is obvious.
}),
normalized in a funny way.  Namely, let
$\nu=Z(\bigcirc)\in\calA$ be the Kontsevich integral of the
unknot\footnote{
  See~\cite{Bar-NatanGaroufalidisRozanskyThurston:WheelsWheeling}
  for the conjectured value of this invariant,
  and~\cite{Bar-NatanLeThurston:InPreparation} for the proof.
},
and let $\Delta_X:\calA\to\calA(\uparrow_X)$ be the ``$X$-cabling''
map that replaces the single directed line in $\calA$ by $n$ directed
lines labeled by the elements of $X$, and sums over all possible ways
of lifting each vertex on the directed line to its $n$ clones (see
e.g.~\cite{Bar-Natan:NAT, LeMurakami:Universal}). Set
\begin{equation} \lbl{eq:ZcheckDef}
  \Zcheck(T)=\nu^{\otimes n}\cdot\Delta_X(\nu)\cdot Z(L)
\end{equation}
for any $X$-marked framed pure tangle $L$, using the action of
$\calA^{\otimes n}$ on $\calA(\uparrow_X)$ defined by sticking
any $n$ diagrams in $\calA$ on the $n$ components of the skeleton of a
diagram in $\calA(\uparrow_X)$. In~\eqref{eq:ZcheckDef} the
factor of $\Delta_X(\nu)$ is easily explained; it accounts for the
difference between the Kontsevich integral of pure tangles and of their
closures. The other factor, $\nu^{\otimes n}$, is the surprising
discovery of~\cite{LeMurakami2Ohtsuki:3Manifold}: it makes $\Zcheck$
better behaved relative to the second Kirby move.
\end{definition}

\begin{definition} \lbl{def:sigma}
The map $\sigma$, first defined in~\cite{Bar-Natan:Homotopy}, is a
simple generalization of the formal PBW map $\sigma:\calA\to\calB$
of~\cite{Bar-Natan:Vassiliev}. It is more easily described through its
inverse $\chi$. If $C\in \calB(X)$ is an $X$-marked uni-trivalent diagram
with $k_x$ legs marked $x$ for any $x\in X$, then
$\chi(C)\in\calA(\uparrow_X)$ is the average\footnote{
  Notice that the normalization is different than
  in~\cite{Bar-Natan:Vassiliev, Bar-Natan:Homotopy} where by careless
  design a sum was used instead of an average.
}
of the $\prod_x k_x!$ ways of attaching the legs of $C$ to $n$ labeled
vertical arrows (labeled by the elements of $X$), attaching legs marked
by $x$ only to the $x$-labeled arrow, for all $x\in X$.
\end{definition}

The new map, $\Gint:\calB(X)\to\calA(\emptyset)$ is only partially
defined.  It's domain is the set $\calB^{FG}(X)$ of ``non-degenerate
perturbed Gaussians'':

\begin{definition} An element $G\in\calB(X)$ is a ``non-degenerate
perturbed Gaussian'' if it is of the form
\[
  G=PG\udot\exp_\udot\left(\frac12\sum_{x,y}l_{xy}\chordn{x}{y}\right)
\]
for some invertible matrix $(l_{xy})$, where
$P:\calB(X)\to\calB^+(X)\subset\calB(X)$ is the natural projection onto the
space $\calB^+(X)$ of uni-trivalent diagrams that have at least one trivalent
vertex on each connected component, and $\chordn{x}{y}$ denotes the
only connected uni-trivalent diagrams that have no trivalent vertices. It is
clear that the matrix $(l_{xy})$ is determined by $G$ if $G$ is of that
form. Call $(l_{xy})$ the ``covariance matrix'' of $G$.
\end{definition}

The last name is explained by claim~\ref{claim:ZisPG}, which the philosophy
majors have stated clearly enough. That claim also says that the following
definition applies to $\sigma\Zcheck(T)\in\calB^{FG}(X)$ for any regular pure
tangle $L$:

\begin{definition} \lbl{def:Gint}
Let $G$ be a non-degenerate perturbed Gaussian with covariance matrix
$(l_{xy})$, and let $(l^{xy})$ be the inverse covariance matrix.  Set
\[ \Gint G=
  \left\langle
    \exp_\udot\left(
      -\frac12\sum_{x,y}l^{xy}\chordu{\partial_x}{\partial_y}
    \right),
  PG
  \right\rangle.
\]
Here the pairing $\langle\cdot,\cdot\rangle:
\calB(\partial_X)\otimes\calB^+(X)\to\calA(\emptyset)$ (where
$\partial_X=\{\partial_x:x\in X\}$ is a set of ``dual'' variables) is
defined by
\[ \langle C_1,C_2\rangle=
  \left(\parbox{3.3in}{
    sum of all ways of gluing the $\partial_x$-marked legs of $C_1$ to the
    $x$-marked legs of $C_2$, for all $x\in X$
  }\right).
\] 
This sum is of course $0$ if the numbers of $x$-marked legs do not
match. If the numbers of legs do match and each diagram has $k_x$ legs
marked by $x$ for $x\in X$, the sum is a sum of $\prod_x k_x!$
terms.
\end{definition}

\subsubsection{Normalization.} \lbl{subsubsec:Normalization} In part II of
this series we will show how the philosophy of section~\ref{subsec:WhatIf}
can be made rigorous, implying the following proposition:

\begin{proposition} \lbl{prop:KirbyII} (Proof
in~\cite[section~\ref{II-sec:Invariance}]{Aarhus:II})
The regular pure tangle invariant $\Aa_0$ descends to an invariant
of regular links and as such it is invariant under the second Kirby move
(figure~\ref{KirbyII}).
\end{proposition}

If we want an invariant of rational homology spheres, we still have to
fix $\Aa_0$ to satisfy the first Kirby move (figure~\ref{KirbyI}). This
is done in a standard way, similar to the way the Kauffman bracket
is tweaked to satisfy the first Reidemeister move. The trick is to
multiply the relatively complicated $\Aa_0$ by a much simpler invariant
of regular links, that has an opposite behavior under the first Kirby
move and is otherwise uninteresting. The result is invariant under both
Kirby moves, and by conservation of interest, it is as interesting as
the original $\Aa_0$:

\begin{definition} Let $U_+$ be the unknot with framing $+1$, and let $U_-$
be the unknot with framing $-1$. Let $L$ be a regular link, and let
$\sigma_+$ ($\sigma_-$) be the number of positive (negative) eigenvalues of
the linking matrix of $L$ (note that $\sigma_\pm$ are invariant under the
second Kirby move, which acts on the linking matrix by a similarity
transformation). Let the \Arhus{} integral $\Aa(L)$ of $L$ be
\[ \Aa(L)=\Aa_0(U_+)^{-\sigma_+}\Aa_0(U_-)^{-\sigma_-}\Aa_0(L), \]
with all products and powers taken using the disjoint union product of
$\calA(\emptyset)$.
\end{definition}

\begin{theorem} (Proof
in~\cite[section~\ref{II-subsec:FirstKirby}]{Aarhus:II})
$\Aa$ is invariant under the first Kirby move as well, and hence it is
an invariant of rational homology spheres.
\end{theorem}

\subsection{The main properties of the \Arhus{} integral} 

\subsubsection{First property of $\Aa$: it has a conceptual foundation} This
property was already proven in section~\ref{sec:intro1}. It is the main
reason why the proofs of all other properties (and of
proposition~\ref{prop:KirbyII}) are relatively simple.

\subsubsection{Second property of $\Aa$: it is universal}
\lbl{subsubsec:Universal}
Like there are finite-type (Vassiliev) invariants of knots (see
e.g.~\cite{Bar-Natan:Vassiliev, Birman:Bulletin, BirmanLin:Vassiliev,
Goussarov:New, Goussarov:nEquivalence, Kontsevich:Vassiliev,
Vassiliev:CohKnot, Vassiliev:Book} and~\cite{Bar-Natan:VasBib}), so there
are finite-type (Ohtsuki) invariants of integer homology spheres (see
e.g.~\cite{Ohtsuki:IntegralHomology, GaroufalidisOhtsuki:3ManifoldsIII,
LeMurakamiOhtsuki:Universal, Le:UniversalIHS} and~\cite{Bar-Natan:VasBib}).
These invariants have a rather simple definition, and just as in the case
of knots, they seem to be rather powerful, though precisely how powerful
they are we still don't know. We argue that the \Arhus{} integral $\Aa$ 
plays in the theory of Ohtsuki invariants the same role as the Kontsevich
integral plays in the theory of Vassiliev invariants. Namely, that it is
a ``universal Ohtsuki invariant''. However, manifold invariants are
somewhat more subtle than knot invariants, and the proper definition of
universality is less transparent (see~\cite{Ohtsuki:IntegralHomology,
GaroufalidisOhtsuki:3ManifoldsIII, Le:UniversalIHS}):

\begin{definition} An invariant $U$ of integer homology spheres with
values in $\calA(\emptyset)$ is a ``universal Ohtsuki invariant'' if
\begin{enumerate}
\item The degree $m$ part $U^{(m)}$ of $U$ is of Ohtsuki type $3m$
   (\cite{Ohtsuki:IntegralHomology}).
\item If $OGL$ denotes the Ohtsuki-Garoufalidis-Le map\footnote{
    Nomenclatorial justification: Ohtsuki~\cite{Ohtsuki:IntegralHomology}
    implicitly considered a map similar to $OGL$, with alternating
    summation over $2^{\{\text{edges}\}}$ rather than over
    $2^{\{\text{vertices}\}}$. Later, Garoufalidis and
    Ohtsuki~\cite{GaroufalidisOhtsuki:3ManifoldsIII} introduced ``white
    vertices'', which amount to an alternating summation over
    $2^{\{\text{edges}\}}\times 2^{\{\text{vertices}\}}$. Finally,
    Le~\cite[lemma~5.1]{Le:UniversalIHS} noticed that in this context
    the alternating summation over $2^{\{\text{edges}\}}$ is superfluous,
    leaving us with the definition presented here.
  },
  defined in figure~\ref{fig:OGLMap}, from manifold diagrams to formal
  linear combinations of unit framed algebraically split links in $S^3$,
  and $S$ denotes the surgery map from such links to integer homology
  spheres, then
  \[ (U\circ S\circ OGL)(D)=D+(\text{higher degree diagrams})\qquad
     (\text{in }\calA(\emptyset))
  \]
  whenever $D$ is a manifold diagram (we implicitly linearly extend $S$
  and $U$, to make this a meaningful equation).
\end{enumerate}
\end{definition}

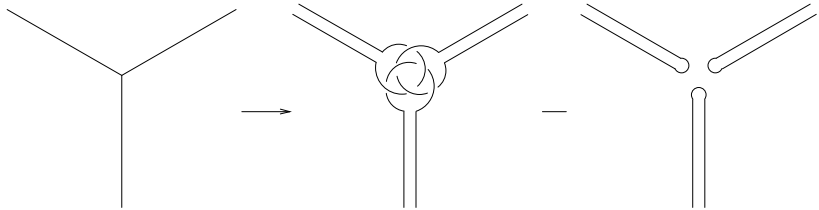
\begin{figure}[htpb]
\[ \printname{OGLMap}
	\setlength{\unitlength}{0.5\standardunitlength}
	\begin{array}{c}  \hspace{-1.7mm}
        	\raisebox{-8pt}{\begingroup\makeatletter\ifx\SetFigFont\undefined
\def\x#1#2#3#4#5#6#7\relax{\def\x{#1#2#3#4#5#6}}%
\expandafter\x\fmtname xxxxxx\relax \def\y{splain}%
\ifx\x\y   
\gdef\SetFigFont#1#2#3{%
  \ifnum #1<17\tiny\else \ifnum #1<20\small\else
  \ifnum #1<24\normalsize\else \ifnum #1<29\large\else
  \ifnum #1<34\Large\else \ifnum #1<41\LARGE\else
     \huge\fi\fi\fi\fi\fi\fi
  \csname #3\endcsname}%
\else
\gdef\SetFigFont#1#2#3{\begingroup
  \count@#1\relax \ifnum 25<\count@\count@25\fi
  \def\x{\endgroup\@setsize\SetFigFont{#2pt}}%
  \expandafter\x
    \csname \romannumeral\the\count@ pt\expandafter\endcsname
    \csname @\romannumeral\the\count@ pt\endcsname
  \csname #3\endcsname}%
\fi
\fi\endgroup
{\renewcommand{\dashlinestretch}{30}
\begin{picture}(10118,2580)(0,-10)
\put(5009.269,1526.423){\arc{664.088}{5.5962}{7.5267}}
\put(5024.929,1538.786){\arc{575.985}{3.3098}{5.0341}}
\put(5003.500,1474.500){\arc{530.330}{1.7127}{2.9997}}
\path(4966,1212)(4966,12)
\path(5116,1212)(5116,12)
\put(5177.159,1701.017){\arc{662.975}{3.5023}{5.4334}}
\put(5159.031,1708.028){\arc{575.878}{1.2162}{2.9387}}
\put(5224.344,1721.431){\arc{531.244}{5.9030}{7.1854}}
\path(5471,1820)(6510,2420)
\path(5396,1950)(6435,2550)
\put(4941.361,1758.389){\arc{661.861}{1.4049}{3.3412}}
\put(4944.438,1739.488){\arc{576.649}{5.4070}{7.1273}}
\put(4900.431,1789.344){\arc{531.244}{3.8102}{5.0926}}
\path(4692,1954)(3653,2553)
\path(4617,1824)(3578,2424)
\put(8641.000,1418.250){\arc{187.500}{2.4981}{6.9267}}
\path(8566,12)(8566,1362)
\path(8716,12)(8716,1362)
\put(8850.516,1784.740){\arc{187.929}{0.4017}{4.8349}}
\path(10106,2423)(8937,1748)
\path(10031,2553)(8862,1878)
\put(8428.484,1782.740){\arc{187.929}{4.5899}{9.0231}}
\path(7248,2551)(8417,1876)
\path(7173,2421)(8342,1746)
\path(1441,12)(1441,1662)
\path(2870,2487)(1441,1662)
\path(12,2487)(1441,1662)
\path(2941,1212)(3541,1212)
\path(3421.000,1182.000)(3541.000,1212.000)(3421.000,1242.000)
\path(6691,1212)(6991,1212)
\end{picture}
} }
        	\hspace{-1.9mm}
	\end{array}
 \]
\caption{
  The $OGL$ map: Take a manifold diagram $D$, embed it in $S^3$ in
  some fixed way of your preference, double every edge, replace every
  vertex by a difference of two local pictures as shown here, and put a
  $+1$ framing on each link component you get. The result is a certain
  alternating sum of $2^v$ links with $e$ components each, where $v$
  and $e$ are the numbers of vertices and edges of $D$, respectively.
}
\lbl{fig:OGLMap}
\end{figure}

\begin{theorem}\label{thm:universal} (Proof
in~\cite[section~\ref{II-subsec:Universal}]{Aarhus:II})
Restricted to integer homology spheres, $\Aa$ is a universal Ohtsuki
invariant.
\end{theorem}

\begin{corollary} \lbl{cor:UniversallityMeans}
\begin{enumerate}
\item $\Aa$ is onto $\calA(\emptyset)$.
\item All Ohtsuki invariants factor through the map $\Aa$.
\item The dual of $\calA(\emptyset)$ is the associated graded of the space
  of Ohtsuki invariants (with degrees divided by $3$; recall
  from~\cite{GaroufalidisOhtsuki:3ManifoldsIII} that the associated
  graded of the space of Ohtsuki invariants vanishes in degrees not
  divisible by $3$.).
\end{enumerate}
\end{corollary}

In view of Le~\cite{Le:UniversalIHS}, this theorem and corollary follow
from the fact (discussed below) that the \Arhus{} integral computes
the $\LMO$ invariant, and from the universality of the $\LMO$
invariant (\cite{Le:UniversalIHS}). But as the definitions of the two
invariants are different, it is nice to have independent proofs of the
main properties.

From corollary~\ref{cor:UniversallityMeans} and the computation of the
low degree parts of $\calA(\emptyset)$ in~\cite{Bar-Natan:Computations}
and~\cite{Kneissler:Twelve} it follows that the low-degree dimensions
of the associated graded of the space of Ohtsuki invariants are given
by the table below. The last row of this table lists the dimensions of
``primitives'' --- multiplicative generators of the algebra of Ohtsuki
invariants.
\begin{center}
\begin{tabular}{|c|r|r|r|r|r|r|r|r|r|r|r|r|}
  \hline
  degree ($3m$) &
    0 & 3 & 6 & 9 & 12 & 15 & 18 & 21 & 24 & 27 & 30 &  33 \\
  \hline
  dimension &
    1 & 1 & 2 & 3 &  6 &  9 & 16 & 25 & 42 & 50 & 90 & 146 \\
  \hline
  dimension of primitives &
    0 & 1 & 1 & 1 &  2 &  2 &  3 &  4 &  5 &  6 &  8 &   9 \\
  \hline
\end{tabular}
\end{center}

\subsubsection{Third property of $\Aa$: it computes the
Le-Murakami-Ohtsuki ($\LMO$) invariant.}

\begin{theorem} \lbl{thm:ArhusIsLMO} (Proof in~\cite{Aarhus:III})
$\Aa$ and the invariant $\widehat{\LMO}$ defined
in~\cite{LeMurakamiOhtsuki:Universal} are equal.\footnote{
  We refer to the invariant $\hat\Omega$
  of~\cite[section~6.2]{LeMurakamiOhtsuki:Universal}. Le, Murakami,
  and Ohtsuki also consider in~\cite{LeMurakamiOhtsuki:Universal} an
  invariant $\Omega$, which we usually denote by $\LMO$. The invariant
  $\LMO$ is defined for general 3-manifolds, and on rational homology
  spheres it differs from $\widehat{\LMO}$ only by a normalization.
}
\end{theorem}

\subsubsection{Fourth property of $\Aa$: it recovers the Rozansky and
Ohtsuki invariants.} \lbl{subsubsec:RelationWithRozansky}

In~\cite{Rozansky:RationalHomology, Rozansky:LinksAndRHS}, Rozansky
shows how to construct a ``perturbative'' invariant of rational
homology spheres, valued in the space of power series in some formal
parameter, which we call $\hbar$. His construction is associated with
the Lie algebra $sl(2)$, but it can be easily extended
(see~\cite{Aarhus:IV}) for other semi-simple Lie algebras $\frakg$. We
call the resulting Rozansky invariant\footnote{
  While the third author participated in writing this paper, he did not
  participate in choosing this terminology.
} $R_\frakg$.

\begin{theorem} \lbl{thm:RecoversRozansky} (Proof in~\cite{Aarhus:IV})
For any semi-simple Lie algebra $\frakg$,
\[ R_\frakg = \Tg \circ \hbar^{\deg} \circ \Aa \]
where $\Tg:\calA(\emptyset)\to\mathbb C$ is the operation of replacing
the vertices and edges of a trivalent graph by the structure constants
of $\frakg$ and the metric of $\frakg$, as in Section~\ref{sec:intro1}
and as in~\cite{Bar-Natan:Vassiliev}, and $\hbar^{\deg}$ is the
operator that multiplies each diagram $D$ in $\calA(\emptyset)$ by
$\hbar$ raised to the degree of $D$.
\end{theorem}

\begin{corollary} \lbl{cor:LMORecoversR}
$\LMO$ recovers $R_\frakg$ for any $\frakg$. In particular,
by~\cite{Rozansky:PrimeK}, $\LMO$ recovers the ``$p$-adic'' invariants
of~\cite{Ohtsuki:PolynomialIntegral,Ohtsuki:PolynomialRational}.
\end{corollary}

The last statement was proven in the case of $\frakg=sl(2)$ by
Ohtsuki~\cite{Ohtsuki:CombinatorialQuantumMethod}.

\section{Frequently Asked Questions}

Let us answer some frequently asked questions. \ifenroute{Most questions
will be answered anyway later in this series, but with the taxi driver
already pushing you out the door, we are sure you would appreciate a
digest. Hold on just a bit more!}{More detailed answers to most of these
questions will be given in the the later parts of this series.}

\begin{question} Your construction uses the Chern-Simons path
integral (at least ideologically). So in what way is your construction
different than the original construction of 3-manifold invariants by
Witten~\cite{Witten:JonesPolynomial}?
\end{question}

\Answer Our path integral is over connections on $S^3$, rather than over
connections on an arbitrary 3-manifold. This means that we can replace
the path integral by any well-behaved universal Vassiliev invariant,
and get a rigorous result.

\begin{question} So what is the relation between your construction and
Witten's?
\end{question}

\Answer One answer is that the \Arhus{} integral should somehow be
related to the $k\to\infty$ asymptotics of the Witten invariants. A more
precise statement is theorem~\ref{thm:RecoversRozansky}; recall that
Rozansky conjectures (and demonstrates in some cases) that his invariant
is related to the trivial connection contribution to the $k\to\infty$
asymptotics of the Witten invariants (in their Reshetikhin-Turaev guise,
see~\cite{Rozansky:RationalHomology, Rozansky:LinksAndRHS}). But perhaps
a more fair answer is {\em we don't know}. There ought to be a direct
path-integral way to see the relation between integration over all
connections on some 3-manifold, and integration over connections on
$S^3$ followed by ``integration of the holonomies'' in the sense of
section~\ref{subsec:WhatIf}. But we don't know this way.

\begin{question} Which one is more general?
\end{question}

\Answer Neither one. Assuming all relevant conjectures, $\Aa$ only
sees the $k\to\infty$ limit, and only ``in the vicinity of the trivial
connection''. But this means it sees a splitting (trivial vs.\ other
flat connections) that the Witten invariants don't see. Also, by
Vogel~\cite{Vogel:Structures}, we know that $\calA(\emptyset)$ ``sees''
more than all semi-simple super Lie algebras see, while the $k\to\infty$
limit of the Witten invariants is practically limited to semi-simple
super Lie algebras. On the other hand, there are Witten-like theories
with finite gauge groups, see e.g.~\cite{FreedQuinn:FiniteGaugeGroup},
which have no parallel in the $\Aa$ world.

\begin{question} What is the relation between the \Arhus{}
integral and the Axelrod-Singer perturbative 3-manifold
invariants~\cite{AxelrodSinger:PerturbationTheoryI,
AxelrodSinger:PerturbationTheoryII} and/or the Kontsevich ``configuration
space integrals''~\cite{Kontsevich:FeynmanDiagrams}?
\end{question}

\Answer We expect the \Arhus{} integral to be the same as Kontsevich's
configuration space integrals and as the formal (no-Lie-algebra) version
of the Axelrod-Singer invariants, perhaps modulo some minor corrections
(in both cases).

\begin{question} What is the relation between the \Arhus{}
integral and the $\LMO$ invariant of Le, Murakami, and
Ohtsuki~\cite{LeMurakamiOhtsuki:Universal}?
\end{question}

\Answer They are the same up to a normalization whenever $\Aa$ is
defined. See theorem~\ref{thm:ArhusIsLMO}.

\begin{question} So what did you add?
\end{question}

\Answer A conceptual construction, and (thus) a better understanding
of the relation between the $\LMO$ invariant/\Arhus{}
integral and the Rozansky and Ohtsuki invariants. See
section~\ref{subsubsec:RelationWithRozansky}.

\begin{question} Does it justify a new name?
\end{question}

\Answer It is rather common in mathematics that different names are used
to describe the same thing, or almost the same thing, depending on the
context or the specific construction. See for example the Kauffman bracket
and the Jones polynomial, the Reidemeister and the Ray-Singer torsions,
and the $\check{\mathrm C}$ech-de-Rham-singular-simplicial cohomology.
Perhaps the name ``\Arhus{} integral'' should only be used when the
construction is explicitly relevant, with the names ``Axelrod-Singer
invariants'', ``Kontsevich's configuration space integrals'', and ``$\LMO$
invariant'' marking the other constructions. When only the functionality
(i.e., theorem~\ref{thm:universal}) matters, the name ``$\LMO$
invariant'' seems most appropriate, as theorem~\ref{thm:universal} was
first considered and proven in the $\LMO$ context.

\begin{question} Didn't Reshetikhin once consider a construction similar
to yours?
\end{question}

\Answer Yes he did, but he never completed his work. Our work was done
independently of his, though in some twisted way it was initiated
by his.  Indeed, it was Reshetikhin's ideas that led one of us
(Rozansky) to study what he called ``the Reshetikhin formula'', and
that led him to discover his ``trivial connection contribution to
the Reshetikhin-Turaev invariants'' (the Rozansky invariants, see
section~\ref{subsubsec:RelationWithRozansky}). The \Arhus{} integral
was discovered by ``reverse engineering'' starting from the Rozansky
invariants --- we first found a diagram-valued invariant that satisfies
theorem~\ref{thm:RecoversRozansky} (working with the ideas we discussed
in~\cite{Bar-NatanGaroufalidisRozanskyThurston:WheelsWheeling}), and
only then we realized that our invariant has the simple interpretation
discussed in section~\ref{sec:intro1}. The result, the \Arhus{} integral,
still carries some affinity to Reshetikhin's construction.

\begin{question} What is the relation between the \Arhus{} integral and the
$p$-adic 3-manifold invariants considered by
Ohtsuki~\cite{Ohtsuki:PolynomialIntegral}?
\end{question}

\Answer See corollary~\ref{cor:LMORecoversR}.

\begin{question} How powerful is $\Aa$?
\end{question}

\Answer It is a ``universal Ohtsuki invariant'' (see
section~\ref{subsubsec:Universal}). In particular, as the Casson invariant
is Ohtsuki-finite-type (see~\cite{Ohtsuki:IntegralHomology}), $\Aa$ is
stronger than the Casson invariant. In terms of the table following
corollary~\ref{cor:UniversallityMeans}, the Casson invariant is just
the degree $3$ primitive, and there are many more. But how powerful
$\Aa$ {\em really} is, how useful it can really be, we don't know. Can
anybody answer that question for the Jones polynomial?

\begin{question} Can you say anything about 3-manifolds with embedded
links?
\end{question}

\Answer Everything works in that case too.
See~\cite[section~\ref{II-subsec:EmbeddedLinks}]{Aarhus:II}.

\begin{question} Can you recommend a good restaurant in \Arhus{}?
\end{question}

\Answer Asian Grill, N{\o}rregade 57, 8000 \Arhus{} C,
tel.~(45)~86~13~86~34. ``The only restaurant in \Arhus{} where the water
isn't vegetarian''.


\section{Appendix: Gaussian integration: a quick refresher.}
\lbl{GaussianDummies}
Whether or not you've seen perturbed Gaussian integration before, you
surely don't want to waste much energy on calculus tricks. Hence we
include here a refresher that you can swallow whole without wasting
any time (and insert all the right factors of $\pm i$, $2^{\pm 1}$,
$\pi^{\pm 1/2}$, etc.\ if you do have some time to spare).

Let $V$ be a vector space and $dv$ a Lebesgue measure on $V$. A perturbed
Gaussian integral is an integral of the form $I_T=\int_V e^T dv$,
where $T$ is a polynomial (or a power series) on $V$, of some specific
form --- $T$ must be a sum $T=\frac12 Q+P$, where $Q$ is a non-degenerate
`big' quadratic, and $P$ is some (possibly) higher degree `perturbation',
which in some sense should be `small' relative to $Q$. If the perturbation
$P$ is missing, then $I_T$ is a simple un-perturbed Gaussian integral.

The first step is to expand the $e^P$ part to a power series,
\[ I_T=\int_V e^T dv=\sum_{m=0}^\infty\frac{1}{m!}
  \int_V P(v)^m e^{Q(v)/2}dv.
\]
Then, we use some Fourier analysis. Recall that the Fourier transform takes
integration to evaluation at $0$, takes multiplication by a polynomial to
multiple differentiation, and takes a Gaussian to another Gaussian, with
the negative-inverse quadratic form. All and all, we find that
\begin{equation} \lbl{eq:NoIntegration}
   I_T\sim \sum_{m=0}^\infty\frac{1}{m!}
    \left.\partial_{P^m}\exp{\left(
      -\frac12 Q^{-1}(v^\star)
    \right)}\right|_{v^\star=0}.
\end{equation}
The notation here means:
\begin{list}{}{
  \setlength{\leftmargin}{50pt}
  \setlength{\labelwidth}{56pt}
  \setlength{\labelsep}{10pt}
}
\item[$\sim$]
  means equality modulo $2$'s, $\pi$'s, $i$'s, and their likes.
\item[$Q^{-1}(v^\star)$]
  is the inverse of $Q(v)$. It is a quadratic form on $V^\star$, and it
  is evaluated on some $v^\star\in V^\star$.
\item[$\partial_{P^m}$]
  is $P^m$ regarded as a differential operator acting on functions on
  $V^\star$. To see how this works, recall that polynomials on $V$ are
  elements of the symmetric algebra $S^\star(V^\star)$ of $V^\star$, and
  they act on $S^\star(V)$, namely on polynomials (and hence functions)
  on $V^\star$, via the standard `contract as much as you can' action
  $S^k(V^\star)\otimes S^l(V)\to S^{l-k}(V)$.
\end{list}

Here are two ways to look at the result, equation~\eqref{eq:NoIntegration}:
\begin{enumerate}
\item We can re-pack the sum as an exponential and get
  \begin{equation} \lbl{eq:CompactForm}
     I_T\sim\left\langle e^P, e^{-Q^{-1}/2}\right\rangle,
  \end{equation}
  where $\left\langle\cdot,\cdot\right\rangle$ denotes the usual pairing
  $S^\star(V^\star)\otimes S^\star(V)\to{\mathbb C}$.
\item We can expand $e^{Q^{-1}/2}$ as a power series and get
  \begin{equation} \lbl{eq:PairingForm}
    I_T\sim\left\langle
      \sum_{m=0}^\infty\frac{P^m}{m!},
      \sum_{n=0}^\infty\frac{(-Q^{-1}/2)^n}{n!}
    \right\rangle.
  \end{equation}
  This last equation has a clear combinatorial interpretation: Take an
  arbitrary number of unordered copies of $P$ and an arbitrary number
  of unordered copies of $Q^{-1}$. If the total degrees happen to be the
  same, you can contract them and get a number. Sum all the numbers you
  thus get, and you've finished computing $I_T$.
\end{enumerate}

For the sake of concreteness, let us play with the example
$T=\frac12 Q+P_2+P_{4,1}+P_{4,2}$, where $Q$ is the big quadratic, $P_2$
is another quadratic which is regarded as a perturbation, and $P_{4,1}$
and $P_{4,2}$ are two additional quartic perturbation terms. It is
natural to represent each of these terms by a picture of an animal
(usually connected) with as many legs as its degree. This is because
$P_{4,1}$ (say) is in $S^4(V^\star)$.  That is, it is a 4-legged animal
which has to be fed with 4 copies of some vector $v$ to produce the
number $P_{4,1}(v)$:
\[ (P_{4,1},v)\mapsto P_{4,1}(v) \qquad\text{is}\qquad
  \left(\printname{Animal}
	\setlength{\unitlength}{0.35\standardunitlength}
	\begin{array}{c}  \hspace{-1.7mm}
        	\raisebox{-8pt}{%
\begingroup\makeatletter\ifx\SetFigFont\undefined
\def\x#1#2#3#4#5#6#7\relax{\def\x{#1#2#3#4#5#6}}%
\expandafter\x\fmtname xxxxxx\relax \def\y{splain}%
\ifx\x\y   
\gdef\SetFigFont#1#2#3{%
  \ifnum #1<17\tiny\else \ifnum #1<20\small\else
  \ifnum #1<24\normalsize\else \ifnum #1<29\large\else
  \ifnum #1<34\Large\else \ifnum #1<41\LARGE\else
     \huge\fi\fi\fi\fi\fi\fi
  \csname #3\endcsname}%
\else
\gdef\SetFigFont#1#2#3{\begingroup
  \count@#1\relax \ifnum 25<\count@\count@25\fi
  \def\x{\endgroup\@setsize\SetFigFont{#2pt}}%
  \expandafter\x
    \csname \romannumeral\the\count@ pt\expandafter\endcsname
    \csname @\romannumeral\the\count@ pt\endcsname
  \csname #3\endcsname}%
\fi
\fi\endgroup
{\renewcommand{\dashlinestretch}{30}
\begin{picture}(1844,1038)(0,-10)
\path(622,912)(1222,912)(1222,312)
	(622,312)(622,912)
\path(622,312)(622,12)
\path(1222,312)(1222,12)
\path(622,912)	(569.605,945.531)
	(523.223,971.721)
	(481.976,990.790)
	(444.985,1002.956)
	(380.254,1007.463)
	(322.000,987.000)

\path(322,987)	(276.753,959.062)
	(235.552,928.390)
	(164.847,857.416)
	(135.123,816.401)
	(109.006,771.224)
	(86.385,721.528)
	(67.150,666.956)
	(51.192,607.151)
	(38.400,541.756)
	(28.666,470.414)
	(24.910,432.401)
	(21.878,392.768)
	(19.555,351.469)
	(17.928,308.460)
	(16.982,263.697)
	(16.704,217.135)
	(17.081,168.728)
	(18.098,118.434)
	(19.743,66.206)
	(22.000,12.000)

\path(1222,912)	(1274.395,945.531)
	(1320.777,971.721)
	(1362.024,990.790)
	(1399.015,1002.956)
	(1463.746,1007.463)
	(1522.000,987.000)

\path(1522,987)	(1567.249,959.063)
	(1608.451,928.393)
	(1679.159,857.422)
	(1708.883,816.408)
	(1735.001,771.231)
	(1757.623,721.536)
	(1776.858,666.964)
	(1792.816,607.159)
	(1805.607,541.763)
	(1815.341,470.421)
	(1819.096,432.408)
	(1822.128,392.774)
	(1824.450,351.474)
	(1826.077,308.465)
	(1827.022,263.701)
	(1827.299,217.138)
	(1826.922,168.731)
	(1825.903,118.435)
	(1824.258,66.207)
	(1822.000,12.000)

\end{picture}
}
 }
        	\hspace{-1.9mm}
	\end{array}
,v\right)\mapsto\printname{FedAnimal}
	\setlength{\unitlength}{0.35\standardunitlength}
	\begin{array}{c}  \hspace{-1.7mm}
        	\raisebox{-8pt}{\begingroup\makeatletter\ifx\SetFigFont\undefined
\def\x#1#2#3#4#5#6#7\relax{\def\x{#1#2#3#4#5#6}}%
\expandafter\x\fmtname xxxxxx\relax \def\y{splain}%
\ifx\x\y   
\gdef\SetFigFont#1#2#3{%
  \ifnum #1<17\tiny\else \ifnum #1<20\small\else
  \ifnum #1<24\normalsize\else \ifnum #1<29\large\else
  \ifnum #1<34\Large\else \ifnum #1<41\LARGE\else
     \huge\fi\fi\fi\fi\fi\fi
  \csname #3\endcsname}%
\else
\gdef\SetFigFont#1#2#3{\begingroup
  \count@#1\relax \ifnum 25<\count@\count@25\fi
  \def\x{\endgroup\@setsize\SetFigFont{#2pt}}%
  \expandafter\x
    \csname \romannumeral\the\count@ pt\expandafter\endcsname
    \csname @\romannumeral\the\count@ pt\endcsname
  \csname #3\endcsname}%
\fi
\fi\endgroup
{\renewcommand{\dashlinestretch}{30}
\begin{picture}(1972,1326)(0,-10)
\path(750,1200)(1350,1200)(1350,600)
	(750,600)(750,1200)
\path(750,600)(750,300)
\path(1350,600)(1350,300)
\path(750,1200)	(697.605,1233.531)
	(651.223,1259.721)
	(609.976,1278.790)
	(572.985,1290.956)
	(508.254,1295.463)
	(450.000,1275.000)

\path(450,1275)	(404.753,1247.062)
	(363.552,1216.390)
	(292.847,1145.416)
	(263.123,1104.401)
	(237.006,1059.224)
	(214.385,1009.528)
	(195.150,954.956)
	(179.192,895.151)
	(166.400,829.756)
	(156.666,758.414)
	(152.910,720.401)
	(149.878,680.768)
	(147.555,639.469)
	(145.928,596.460)
	(144.982,551.697)
	(144.704,505.135)
	(145.081,456.728)
	(146.098,406.434)
	(147.743,354.206)
	(150.000,300.000)

\path(1350,1200)	(1402.395,1233.531)
	(1448.777,1259.721)
	(1490.024,1278.790)
	(1527.015,1290.956)
	(1591.746,1295.463)
	(1650.000,1275.000)

\path(1650,1275)	(1695.249,1247.063)
	(1736.451,1216.393)
	(1807.159,1145.422)
	(1836.883,1104.408)
	(1863.001,1059.231)
	(1885.623,1009.536)
	(1904.858,954.964)
	(1920.816,895.159)
	(1933.607,829.763)
	(1943.341,758.421)
	(1947.096,720.408)
	(1950.128,680.774)
	(1952.450,639.474)
	(1954.077,596.465)
	(1955.022,551.701)
	(1955.299,505.138)
	(1954.922,456.731)
	(1953.903,406.435)
	(1952.258,354.207)
	(1950.000,300.000)

\put(600,0){\makebox(0,0)[lb]{\smash{{\mathmode{v}}}}}
\put(1200,0){\makebox(0,0)[lb]{\smash{{\mathmode{v}}}}}
\put(1800,0){\makebox(0,0)[lb]{\smash{{\mathmode{v}}}}}
\put(0,0){\makebox(0,0)[lb]{\smash{{\mathmode{v}}}}}
\end{picture}
} }
        	\hspace{-1.9mm}
	\end{array}
.
\]

With this in mind, $T$ is represented by a sum of such connected animals:
\[ \printname{Vapor}
	\setlength{\unitlength}{0.35\standardunitlength}
	\begin{array}{c}  \hspace{-1.7mm}
        	\raisebox{-8pt}{\input draws/Vapor.tex }
        	\hspace{-1.9mm}
	\end{array}
. \]
\parbox{3.5in}{
  Exponentiation is done using power series. Each term in $\exp T$
  is some power of $T$ and hence an element of some symmetric power
  $S^k(V^\star)$. As such, it is represented by some sum of pictures, each
  of which is some disjoint union of the connected animals making $T$.
  Namely, $\exp T$ is some sum of `clouds'\ifenroute{ (Dylan thinks
  these are ``swarms of insects with unusual mating practices'')}\ like
  the one on the right.
}\hfill$\printname{Cloud}
	\setlength{\unitlength}{0.35\standardunitlength}
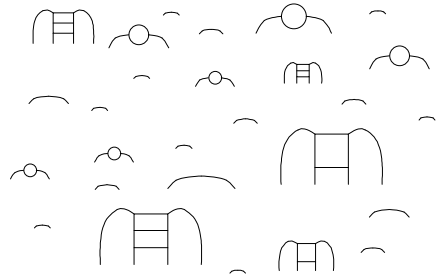
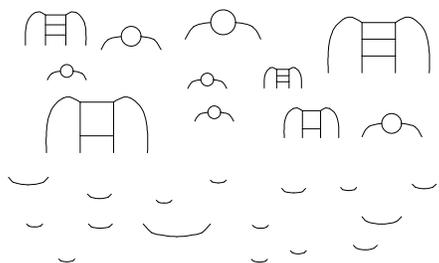
	\begin{array}{c}  \hspace{-1.7mm}
        	\raisebox{-8pt}{\input draws/Cloud.tex }
        	\hspace{-1.9mm}
	\end{array}
$\hfill\ 

\par\noindent
\vskip 1mm
\ \hfill$\printname{SplitCloud}
	\setlength{\unitlength}{0.35\standardunitlength}
	\begin{array}{c}  \hspace{-1.7mm}
        	\raisebox{-8pt}{\input draws/SplitCloud.tex }
        	\hspace{-1.9mm}
	\end{array}
$\hfill\parbox{3.5in}{
  The next step, as seen from equation~\eqref{eq:PairingForm}, is to
  separate $Q$ from the rest, negate it, and invert it. If we think of
  legs pointing down as legs in $V^\star$ and legs pointing up as legs
  in $V$, the result is a sum of `split coulds' like the one on the left.
}
\vskip 1mm

The final step according to equation~\eqref{eq:PairingForm} is to contract
the $(-Q^{-1})$'s with the $P$'s, whenever the degrees allow that. In
pictures, we just connect the down-pointing legs to the up-pointing legs in
all possible ways, and the result is a big sum of diagrams that look like
this:
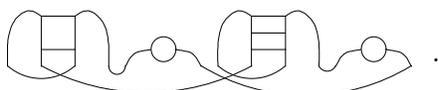
\[ \printname{GluedCloud}
	\setlength{\unitlength}{0.35\standardunitlength}
	\begin{array}{c}  \hspace{-1.7mm}
        	\raisebox{-8pt}{\begingroup\makeatletter\ifx\SetFigFont\undefined
\def\x#1#2#3#4#5#6#7\relax{\def\x{#1#2#3#4#5#6}}%
\expandafter\x\fmtname xxxxxx\relax \def\y{splain}%
\ifx\x\y   
\gdef\SetFigFont#1#2#3{%
  \ifnum #1<17\tiny\else \ifnum #1<20\small\else
  \ifnum #1<24\normalsize\else \ifnum #1<29\large\else
  \ifnum #1<34\Large\else \ifnum #1<41\LARGE\else
     \huge\fi\fi\fi\fi\fi\fi
  \csname #3\endcsname}%
\else
\gdef\SetFigFont#1#2#3{\begingroup
  \count@#1\relax \ifnum 25<\count@\count@25\fi
  \def\x{\endgroup\@setsize\SetFigFont{#2pt}}%
  \expandafter\x
    \csname \romannumeral\the\count@ pt\expandafter\endcsname
    \csname @\romannumeral\the\count@ pt\endcsname
  \csname #3\endcsname}%
\fi
\fi\endgroup
{\renewcommand{\dashlinestretch}{30}
\begin{picture}(7230,1562)(0,-10)
\path(622,1436)(1222,1436)(1222,836)
	(622,836)(622,1436)
\path(622,836)(622,536)
\path(1222,836)(1222,536)
\path(622,1436)	(569.605,1469.531)
	(523.223,1495.721)
	(481.976,1514.790)
	(444.985,1526.956)
	(380.254,1531.463)
	(322.000,1511.000)

\path(322,1511)	(276.753,1483.062)
	(235.552,1452.390)
	(164.847,1381.416)
	(135.123,1340.401)
	(109.006,1295.224)
	(86.385,1245.528)
	(67.150,1190.956)
	(51.192,1131.151)
	(38.400,1065.756)
	(28.666,994.414)
	(24.910,956.401)
	(21.878,916.768)
	(19.555,875.469)
	(17.928,832.460)
	(16.982,787.697)
	(16.704,741.135)
	(17.081,692.728)
	(18.098,642.434)
	(19.743,590.206)
	(22.000,536.000)

\path(1222,1436)	(1274.395,1469.531)
	(1320.777,1495.721)
	(1362.024,1514.790)
	(1399.015,1526.956)
	(1463.746,1531.463)
	(1522.000,1511.000)

\path(1522,1511)	(1567.249,1483.063)
	(1608.451,1452.393)
	(1679.159,1381.422)
	(1708.883,1340.408)
	(1735.001,1295.231)
	(1757.623,1245.536)
	(1776.858,1190.964)
	(1792.816,1131.159)
	(1805.607,1065.763)
	(1815.341,994.421)
	(1819.096,956.408)
	(1822.128,916.774)
	(1824.450,875.474)
	(1826.077,832.465)
	(1827.022,787.701)
	(1827.299,741.138)
	(1826.922,692.731)
	(1825.903,642.435)
	(1824.258,590.207)
	(1822.000,536.000)

\put(2797,836){\ellipse{450}{450}}
\path(2572,836)	(2510.982,828.237)
	(2458.441,820.452)
	(2413.498,812.427)
	(2375.275,803.941)
	(2315.472,784.710)
	(2272.000,761.000)

\path(2272,761)	(2206.262,688.576)
	(2168.205,624.944)
	(2146.340,583.965)
	(2122.000,536.000)

\path(3022,836)	(3083.018,828.237)
	(3135.559,820.452)
	(3180.502,812.427)
	(3218.725,803.941)
	(3278.528,784.710)
	(3322.000,761.000)

\path(3322,761)	(3387.737,688.576)
	(3425.795,624.944)
	(3447.660,583.965)
	(3472.000,536.000)

\put(6547,836){\ellipse{450}{450}}
\path(6322,836)	(6260.982,828.237)
	(6208.441,820.452)
	(6163.498,812.427)
	(6125.275,803.941)
	(6065.472,784.710)
	(6022.000,761.000)

\path(6022,761)	(5956.263,688.576)
	(5918.205,624.944)
	(5896.340,583.965)
	(5872.000,536.000)

\path(6772,836)	(6833.018,828.237)
	(6885.559,820.452)
	(6930.502,812.427)
	(6968.725,803.941)
	(7028.528,784.710)
	(7072.000,761.000)

\path(7072,761)	(7137.737,688.576)
	(7175.795,624.944)
	(7197.660,583.965)
	(7222.000,536.000)

\path(4373,1435)(4973,1435)(4973,835)
	(4373,835)(4373,1435)
\path(4373,835)(4373,535)
\path(4973,835)(4973,535)
\path(4372,1136)(4972,1136)
\path(4373,1435)	(4320.605,1468.531)
	(4274.223,1494.721)
	(4232.976,1513.790)
	(4195.985,1525.956)
	(4131.254,1530.463)
	(4073.000,1510.000)

\path(4073,1510)	(4027.753,1482.062)
	(3986.552,1451.390)
	(3915.847,1380.416)
	(3886.123,1339.401)
	(3860.006,1294.224)
	(3837.385,1244.528)
	(3818.150,1189.956)
	(3802.192,1130.151)
	(3789.400,1064.756)
	(3779.666,993.414)
	(3775.910,955.401)
	(3772.878,915.768)
	(3770.555,874.469)
	(3768.928,831.460)
	(3767.982,786.697)
	(3767.704,740.135)
	(3768.081,691.728)
	(3769.098,641.434)
	(3770.743,589.206)
	(3773.000,535.000)

\path(4973,1435)	(5025.395,1468.531)
	(5071.777,1494.721)
	(5113.024,1513.790)
	(5150.015,1525.956)
	(5214.746,1530.463)
	(5273.000,1510.000)

\path(5273,1510)	(5318.249,1482.063)
	(5359.451,1451.393)
	(5430.159,1380.422)
	(5459.883,1339.408)
	(5486.001,1294.231)
	(5508.623,1244.536)
	(5527.858,1189.964)
	(5543.816,1130.159)
	(5556.607,1064.763)
	(5566.341,993.421)
	(5570.096,955.408)
	(5573.128,915.774)
	(5575.450,874.474)
	(5577.077,831.465)
	(5578.022,786.701)
	(5578.299,740.138)
	(5577.922,691.731)
	(5576.903,641.435)
	(5575.258,589.207)
	(5573.000,535.000)

\put(621.000,1224.500){\arc{1825.000}{0.8533}{2.2883}}
\put(1971.000,537.000){\arc{300.000}{6.2832}{9.4248}}
\put(5721.000,612.000){\arc{335.410}{0.4636}{2.6779}}
\put(5346.000,3601.286){\arc{7184.838}{1.0217}{2.1199}}
\put(2496.000,3993.250){\arc{7864.169}{1.0737}{2.0679}}
\put(4371.000,1224.500){\arc{1825.000}{0.8533}{2.2883}}
\end{picture}
} }
        	\hspace{-1.9mm}
	\end{array}
. \]
Notice that in these diagrams (commonly referred to as ``Feynman
diagrams'') there are no `free legs' left. Therefore they represent
complete contractions, that is, scalars.

\par\noindent
{\bf Acknowledgement: } The seeds leading to this work were planted when
the four of us (as well as Le, Murakami (H\&J), Ohtsuki, and many other
like-minded people) were visiting \Arhus, Denmark, for a special semester
on geometry and physics, in August 1995. We wish to thank the organizers,
J.~Dupont, H.~Pedersen, A.~Swann and especially J.~Andersen for their
hospitality and for the stimulating atmosphere they created. We also
wish to thank N.~Habegger, M.~Hutchings, T.~Q.~T.~Le, A.~Referee and
N.~Reshetikhin for additional remarks and suggestions, and the Center
for Discrete Mathematics and Theoretical Computer Science at the Hebrew
University for financial support.

\ifx\undefined\bysame
        \newcommand{\bysame}{\leavevmode\hbox to3em{\hrulefill}\,}
\fi

\end{document}
\endinput